\documentclass{article}
\usepackage{amsfonts,color}
\usepackage{amsmath}
 \usepackage{epsfig}
\usepackage{lscape}
\usepackage{amssymb}
\usepackage{bbm}
\renewcommand{\epsilon}{\varepsilon}
\newcommand{\R}{\mathbb{R}}

\DeclareMathOperator{\polar}{polar}
\newcommand{\id}{{\boldsymbol{\mathbbm{1}}}}

\allowdisplaybreaks[1]

\makeindex

\def\barr{\begin{array}}

\def\earr{\end{array}}
\def\bec#1{\begin{equation}\label{#1}}
\def\becn{\begin{equation*}}
\def\endec{\end{equation}}
\def\endecn{\end{equation*}}
\def\dd{\displaystyle}

 \usepackage{footnote}

\makeatletter
\let\@fnsymbol\@arabic
\makeatother

\begin{document}
\title{Refined dimensional reduction for isotropic elastic Cosserat  shells with initial curvature}
\author{ Mircea B\^irsan\thanks{Mircea B\^irsan, \ \  Lehrstuhl f\"{u}r Nichtlineare Analysis und Modellierung, Fakult\"{a}t f\"{u}r Mathematik,
Universit\"{a}t Duisburg-Essen, Thea-Leymann Str. 9, 45127 Essen, Germany; and  Alexandru Ioan Cuza University of Ia\c si, Department of Mathematics,  Blvd.
Carol I, no. 11, 700506 Ia\c si,
Romania;  email: mircea.birsan@uni-due.de}  
\quad 
and  \quad  Ionel-Dumitrel Ghiba\thanks{Ionel-Dumitrel Ghiba, \ \   Alexandru Ioan Cuza University of Ia\c si, Department of Mathematics,  Blvd.
Carol I, no. 11, 700506 Ia\c si,
Romania; and  Octav Mayer Institute of Mathematics of the
Romanian Academy, Ia\c si Branch,  700505 Ia\c si, email:  dumitrel.ghiba@uaic.ro}
 \quad and\quad Robert J. Martin\,\thanks{Robert J. Martin,  \ \ Lehrstuhl f\"{u}r Nichtlineare Analysis und Modellierung, Fakult\"{a}t f\"{u}r
 	Mathematik, Universit\"{a}t Duisburg-Essen,  Thea-Leymann Str. 9, 45127 Essen, Germany, email: robert.martin@uni-due.de}
  \\ and\quad Patrizio Neff\,\thanks{Patrizio Neff,  \ \ Head of Lehrstuhl f\"{u}r Nichtlineare Analysis und Modellierung, Fakult\"{a}t f\"{u}r
  	Mathematik, Universit\"{a}t Duisburg-Essen,  Thea-Leymann Str. 9, 45127 Essen, Germany, email: patrizio.neff@uni-due.de}
}

\maketitle

\begin{center}
\thanks{\textit{Dedicated  to  Sanda Cleja-\c Tigoiu on the occasion of her 70th birthday}}
\end{center}

\begin{abstract}
Using a geometrically motivated 8-parameter ansatz through the thickness, we reduce a three-dimensional shell-like geometrically nonlinear Cosserat material to a fully two-dimensional shell model.
Curvature effects are fully taken into account. For elastic isotropic Cosserat materials, the integration through the thickness can be performed analytically   and a generalized plane stress condition allows for a closed-form expression of the thickness stretch and the nonsymmetric shift of the midsurface in bending. We obtain an explicit  form of the elastic strain energy density for Cosserat shells, including  terms up to order $ O(h^5) $ in the shell thickness $ h $. This energy density is expressed as a quadratic function of the nonlinear elastic shell strain tensor and the bending-curvature tensor, with coefficients depending on the initial curvature of the shell.
\end{abstract}

\section{Introduction}

Nonlinear elastic shell theory is a notoriously difficult subject from the perspective of modelling, analysis and numerical implementation. Some successful  research has been devoted to generalizations of the Reissner-Mindlin kinematics, well known from linear elastic plate models \cite{Neff_Hong_Reissner08}. In his habilitation thesis, the last author of the present article began the modelling and analysis of so-called nonlinear Cosserat shell models, in which a full triad of orthogonal directors, independent of the normal of the shell, is taken into account \cite{Neff_Habil04,Neff_plate04_cmt,Neff_plate07_m3as,Neff_Chelminski_ifb07}. As such, these models fall into the class of 6-parameter shell models, proposed originally by Reissner \cite{Reissner74} and presented in the book of Libai and Simmonds \cite{Libai98} as well as in the works of Pietraszkiewicz and coauthors \cite{Pietraszkiewicz04,Pietraszkiewicz-book04,Eremeyev06}, see also \cite{Tambaca-16,Tambaca-19}.
They are also the preferred models from an engineering point of view, since the independent rotation field allows for transparent coupling between shell and beam parts.

The results of \cite{Neff_Habil04,Neff_plate04_cmt,Neff_plate07_m3as} have been obtained by an 8-parameter ansatz of the deformation through the thickness and consistent analytic integration over the thickness in the case of a flat undeformed shell reference configuration. The mathematical approach developed there also allowed for the first existence proof of minimizers \cite{Neff_plate04_cmt,Neff_plate07_m3as,Birsan-Neff-MMS-2014}. In this paper, we extend the modelling from flat shells to initially curved shells. With appropriate changes, taking into account the geometry of the shell, it is possible to use the same derivation ideas. Our ansatz allows for a consistent shell model up to order $O(h^5)$ in the shell thickness. Interestingly, the contributions of order $h^5$ in the shell energy all depend on the initial curvature of the shell and vanish for a flat-shell. In this case we recover the previously mentioned flat Cosserat shell model \cite{Neff_plate04_cmt}. The $h^5$-contribution does not come with a definite sign, such that the additional terms can be stabilizing as well as destabilizing, depending on the local shell geometry. 

However, all occurring material coefficients of the shell model are uniquely determined from the isotropic three-dimensional Cosserat model and the given initial geometry of the shell. Thus, we fill a certain gap in the general 6-parameter shell theory, which leaves the precise structure of the constitutive equations wide open.

In a future contribution we will try to prove the existence of minimizers for our new Cosserat shell model along the methods outlined in 
\cite{Birsan-Neff-MMS-2014}.

\setcounter{equation}{0}

\section{The three-dimensional Cosserat model in curvilinear coordinates}

Let $\Omega_\xi\subset\mathbb{R}^3$ be the reference configuration of an elastic Cosserat body. The elastic material occupying the domain $\Omega_\xi$ is assumed to be homogeneous and isotropic. A generic point of $\Omega_\xi$ will be denoted by $(\xi_1,\xi_2,\xi_3)$.
The deformation of the body is described by a vectorial map $\boldsymbol{\varphi}_\xi$ (called \emph{deformation}) and  a \textit{microrotation tensor} $\boldsymbol{R}_\xi$:
\begin{equation}\label{e1}
\boldsymbol \varphi_\xi:\Omega_\xi \rightarrow\mathbb{R}^3, \qquad \boldsymbol{R}_\xi:\Omega_\xi \rightarrow \mathrm{SO}(3)\, .
\end{equation}
We denote by $\Omega_c:=\boldsymbol \varphi_\xi(\Omega_\xi)$ the current (deformed) configuration.

Throughout the paper, we adhere to some common notational conventions: boldface letters denote vectors and tensors; the Latin indices $i,j,k,...$ range over the set $\{1,2,3\}$, while the Greek indices $\alpha,\beta,\gamma,...$ take the values $\{1,2\}$. We also employ the Einstein summation convention over repeated indices.

\begin{figure}
	\begin{center}
		\includegraphics[scale=1]{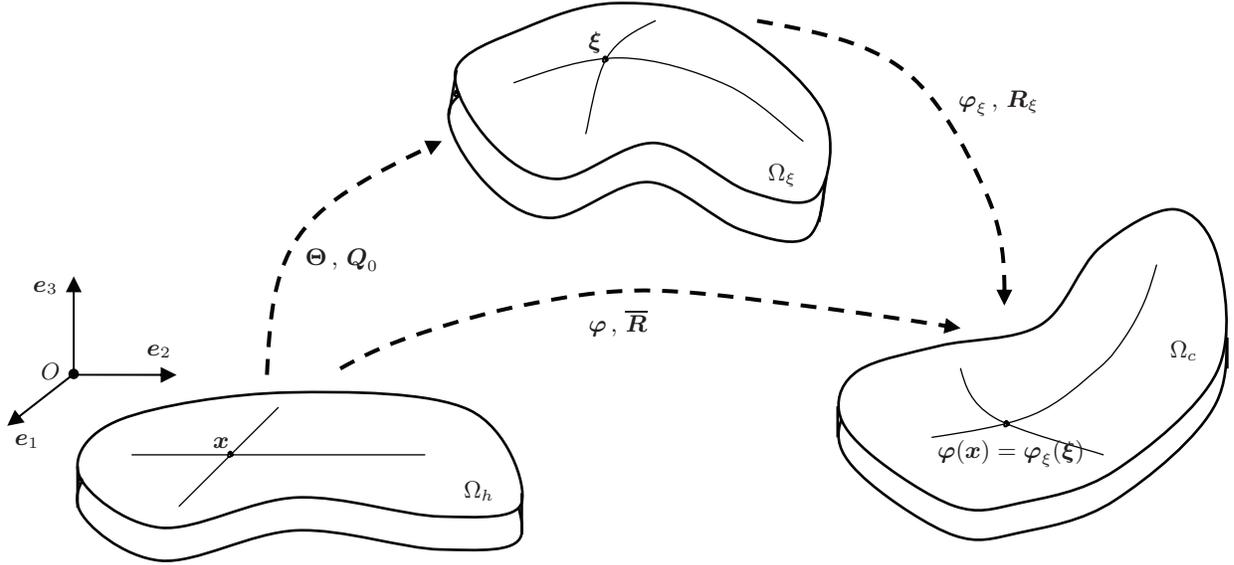}
		\put(-452,72){\small{$O$}} \put(-462,47){\small{$\boldsymbol{e}_1$}}
		\put(-412,81){\small{$\boldsymbol{e}_2$}}
		\put(-455,105){\small{$\boldsymbol{e}_3$}}
		\put(-387,46){\small{$\boldsymbol{x}$}} 
		\put(-292,27){\small{$\Omega_h$}} 
		\put(-177,147){\small{$\Omega_\xi$}} 
		\put(-245,197){\small{$\boldsymbol{\xi}$}} 
		\put(-25,80){\small{$\Omega_c$}} 
		\put(-113,42){\small{$\boldsymbol \varphi(\boldsymbol{x})=\boldsymbol \varphi_\xi(\boldsymbol{\xi})$}} 
		\put(-352,115){\small{$\boldsymbol\Theta \,, \,\boldsymbol Q_0$}} 
		\put(-245,90){\small{$\boldsymbol\varphi \,, \,\overline{\boldsymbol R}$}}
		\put(-105,175){\small{$\boldsymbol\varphi_\xi \,, \,\boldsymbol R_\xi$}} 
		\caption{The reference configuration $\Omega_\xi$ of the shell, the current configuration $\Omega_c$ and the ``fictitious flat Cartesian configuration'' $\Omega_h\,$.}
		\label{Fig1}       % Give a unique label
	\end{center}
\end{figure}

Consider a  parametric representation
\begin{equation}\label{e2}
\boldsymbol\Theta:\Omega_h \rightarrow\Omega_\xi\,, \qquad \boldsymbol\Theta(x_1,x_2,x_3)=(\xi_1,\xi_2,\xi_3)=\xi_i \boldsymbol e_i\,,
\end{equation}
of the domain  $\Omega_\xi$, where $\boldsymbol \Theta$  is a $C^1$-diffeomorphism and $\Omega_h\subset \mathbb{R}^{3}$ is the parameter domain. Here, the domains $\Omega_h $ and $ \Omega_\xi $ are equipped with a right Cartesian coordinate frame with unit vectors $
\boldsymbol e_i$ along the axes $Ox_i\,$, see Figure \ref{Fig1}. Thus, the domain $\Omega_h $ can be viewed as a ``fictitious flat Cartesian configuration'' of the body. 
In view of \eqref{e2}, we can regard $(x_1,x_2,x_3)$ as curvilinear coordinates on $\Omega_\xi\,$. The special form of the parameter domain $\Omega_h$ and of the parametrization function $\boldsymbol \Theta$ appropriate to thin shells will be introduced in the next section.

With respect to these curvilinear coordinates $x_i$ $(i=1,2,3)$, we use the covariant base vectors $\boldsymbol g_i$ and the contravariant base vectors $\boldsymbol g^i$ given by
\begin{equation}\label{e3}
    \boldsymbol g_i=\dfrac{\partial\boldsymbol\Theta }{\partial x_i}\,= \dfrac{\partial \xi_j  }{\partial x_i}
    \,\boldsymbol e_j\,,\qquad 
    \boldsymbol g^i=\nabla_\xi\, x_i\,= \dfrac{\partial x_i }{\partial \xi_j }
    \,\boldsymbol e_j\,,
    \qquad
    \boldsymbol g^j\cdot \boldsymbol g_i=\delta^j_i\,,
\end{equation}
where we denote by $ \cdot $ the scalar product and $\delta^j_i$ is the Kronecker symbol.
Using the parametrization $\boldsymbol \Theta$, the \textit{deformation function}  $\boldsymbol \varphi$ may be written  as  the composition
\begin{equation}\label{e4}
    \boldsymbol \varphi:=\boldsymbol \varphi_\xi\circ\boldsymbol\Theta :\Omega_h \rightarrow\Omega_c\,,\qquad \boldsymbol \varphi(x_1,x_2,x_3): = \boldsymbol \varphi_\xi\big( \boldsymbol\Theta(x_1,x_2,x_3)\big).
\end{equation}
The gradient of $\boldsymbol\Theta$ can be represented as
\begin{equation}\label{e5}
    \nabla_x \boldsymbol\Theta=\dfrac{\partial\boldsymbol\Theta }{\partial x_i}\,\,\otimes \boldsymbol e_i=\boldsymbol g_i\otimes \boldsymbol e_i\qquad\textrm{and} \qquad  \big(\nabla_x \boldsymbol\Theta\big)^{-1} =\boldsymbol e_i\otimes \boldsymbol g^i,
\end{equation}
where $\otimes$ denotes the diadic product. 
The polar decomposition of $\nabla_x \boldsymbol \Theta$ is denoted by
\begin{equation}\label{e6}
\nabla_x \boldsymbol \Theta={\boldsymbol Q}_0 \,\boldsymbol U_0=\polar{(\nabla_x \boldsymbol \Theta)}\,{\boldsymbol U}_0\, ,
\end{equation}
where $\,{\boldsymbol Q}_0=\polar{(\nabla_x \boldsymbol \Theta)}\in \rm{SO}(3 )$ is a pure rotation and $\,{\boldsymbol U}_0\,$ is a positive definite symmetric tensor \cite{Neff-Fischle-14}. Based on the polar factor ${\boldsymbol Q}_0$, we define the orthonormal directors $\boldsymbol d_1^0\,,\,\boldsymbol d_2^0\,,\,\boldsymbol d_3^0$ in the reference configuration $\Omega_\xi$ by
\begin{equation}\label{e7}
    \boldsymbol d_i^0:=  \boldsymbol Q_0\,\boldsymbol e_i\,, \qquad\mbox{i.e.} \qquad   \boldsymbol Q_0=\boldsymbol d_i^0\otimes\boldsymbol e_i\,.
\end{equation}
Furthermore, we introduce the \emph{elastic microrotation} $\boldsymbol Q_e$ as the composition
\begin{equation}\label{e7,5}
     \boldsymbol Q_e= \boldsymbol R_\xi\circ\boldsymbol\Theta :\Omega_h \rightarrow \mathrm{SO}(3 ),\qquad \boldsymbol Q_e (\boldsymbol x):= \boldsymbol R_\xi\big(\boldsymbol\Theta(\boldsymbol x)\big),
\end{equation}
where $\boldsymbol x:=(x_1,x_2,x_3)=x_i \boldsymbol e_i\,$.
Then the orthonormal directors $\boldsymbol d_1\,,\,\boldsymbol d_2\,,\,\boldsymbol d_3$ in the current configuration $\Omega_c$ are given by
\begin{equation}\label{e8}
    \boldsymbol d_i:=  \boldsymbol Q_e\boldsymbol d^0_i \,, \qquad\mbox{i.e.} \qquad    \boldsymbol Q_e=\boldsymbol d_i\otimes\boldsymbol d_i^0\,;
\end{equation}
these directors describe the orientation of the material points of the Cosserat medium and characterize the microrotation field. Moreover, the \emph{total microrotation} $\overline{\boldsymbol R}$ is given by
\begin{equation}\label{e9}
    \overline{\boldsymbol R}:\Omega_h \rightarrow \mathrm{SO}(3 ),\qquad \overline{\boldsymbol R}(\boldsymbol x):=\boldsymbol Q_e(\boldsymbol x)\,\boldsymbol Q_0(\boldsymbol x)\,.
\end{equation}
In particular,
\begin{equation}\label{e10}
    \overline{\boldsymbol R}(\boldsymbol x)= (\boldsymbol d_i\otimes\boldsymbol d_i^0)\, (\boldsymbol d_j^0\otimes\boldsymbol e_j)= \boldsymbol d_i(\boldsymbol x)\otimes\boldsymbol e_i\,.
\end{equation}

\textbf{Remark:}  We write the microrotation  tensor $\overline{\boldsymbol{R}}$  with a superposed  bar in order to distinguish it from the orthogonal factor  $ {\boldsymbol{R}}$ in the classical polar decomposition of the deformation gradient $\boldsymbol{F}=\boldsymbol{R}\,\boldsymbol{U}\,$, which is a common notational convention in nonlinear elasticity.\hfill$ \Box $ \medskip

The \emph{deformation gradient} corresponding to the deformation $\boldsymbol\varphi_\xi$ is given by
\begin{equation}\label{e11}
\boldsymbol F_\xi:=\nabla_\xi\,  \boldsymbol \varphi_\xi= \dfrac{\partial\boldsymbol\varphi_\xi }{\partial \xi_i}\,\,\otimes \boldsymbol e_i\,,
\end{equation}
whereas $ \boldsymbol F $ denotes the gradient of the deformation function $ \boldsymbol \varphi $ with respect to $ \boldsymbol x $, i.e.
\begin{equation}\label{e11,1}
\boldsymbol F:=\nabla_x\,  \boldsymbol \varphi = \dfrac{\partial\boldsymbol\varphi }{\partial x_i}\,\,\otimes \boldsymbol e_i\,.
\end{equation}
Thus
\begin{equation}\label{e11,2}
\boldsymbol F_\xi=\boldsymbol F \big(\nabla_x \boldsymbol\Theta\big)^{-1}
\end{equation}
since, due to \eqref{e3}, \eqref{e5} and the chain rule,
\[  
\boldsymbol F \,= \,\dfrac{\partial\boldsymbol\varphi_\xi }{\partial \xi_j}\,\dfrac{\partial  \xi_j }{\partial x_i}\,\otimes \boldsymbol e_i
\,=\, \dfrac{\partial\boldsymbol\varphi_\xi }{\partial \xi_j}\,(\boldsymbol e_j\cdot \boldsymbol g_i)\,\otimes \boldsymbol e_i 
\,=  \Big(\dfrac{\partial\boldsymbol\varphi_\xi }{\partial \xi_j}\otimes\boldsymbol e_j\Big)( \boldsymbol g_i\otimes \boldsymbol e_i)
= \boldsymbol F_\xi \big(\nabla_x \boldsymbol\Theta\big)
.
\]
By virtue of \eqref{e5} and \eqref{e11,1}, the relation \eqref{e11,2} can also be written as
\begin{equation}\label{e11,25}
\boldsymbol F_\xi =  \Big(\dfrac{\partial\boldsymbol\varphi }{\partial x_i}\otimes\boldsymbol e_i\Big)( \boldsymbol e_j\otimes \boldsymbol g^j)
=  \dfrac{\partial\boldsymbol\varphi }{\partial x_i}\,\,\otimes \boldsymbol g^i =\boldsymbol \varphi,_i\otimes\, \boldsymbol g^i\,,
\end{equation}
where an index $i$ preceded by a comma denotes partial differentiation with respect to $x_i$\,.

For the \emph{non-symmetric Biot-type stretch tensor}
\begin{equation}\label{e11,3}
\overline{\boldsymbol{U}}_\xi:=\boldsymbol R_\xi^T\, \boldsymbol F_\xi \qquad\textrm{and}\qquad 
\overline{\boldsymbol{U}}_e (\boldsymbol x):=  \overline{\boldsymbol{U}}_\xi\big(\boldsymbol\Theta(\boldsymbol x)\big)\,,
\end{equation}
we find, using \eqref{e7,5}, \eqref{e8} and \eqref{e11,2},
\begin{equation}\label{e12}
    \overline{\boldsymbol{U}}_e=\boldsymbol Q_e^T\, \boldsymbol F \big(\nabla_x \boldsymbol\Theta\big)^{-1} = \big(\boldsymbol d_i^0\otimes\boldsymbol d_i\big)\,\big(\boldsymbol \varphi,_j\otimes \boldsymbol g^j \big)= \big(\boldsymbol \varphi,_j\cdot\, \boldsymbol d_i\big)\,
    \big(\boldsymbol d_i^0\otimes \boldsymbol g^j \big).
\end{equation}
We denote by $\id_3= \boldsymbol g_j\otimes \boldsymbol g^j$ the three-dimensional unit tensor and introduce the non-symmetric \emph{strain tensor}
\begin{equation}\label{e13}
    \overline{\boldsymbol{E}}:=\overline{\boldsymbol{U}}_e-\id_3=   \big(\boldsymbol \varphi,_j\cdot\, \boldsymbol d_i- \boldsymbol g_j\cdot \boldsymbol d_i^0\big)\,
    \big(\boldsymbol d_i^0\otimes \boldsymbol g^j \big),
\end{equation}
which is a Lagrangian strain measure for stretch \cite{Pietraszkiewicz09}.
As a Lagrangian strain measure for curvature (orientation change), we employ the so-called \emph{wryness tensor} (see e.g., \cite{Neff_curl08,Birsan-Neff-L58-2017})
\begin{equation}\label{e14}
\boldsymbol \Gamma_\xi:= \mathrm{axl}\big(\boldsymbol R_\xi^T\, \partial_{\xi_i}\boldsymbol R_\xi\big)\otimes \boldsymbol e_i\,,
\end{equation}
where $ \mathrm{axl}(\cdot) $ is the axial vector of any skew-symmetric tensor and
 $ \partial_{\xi_i} $ denotes the partial differentiation with respect to $ \xi_i\, $.
In the following, we will show that
\begin{equation}\label{e14,1}
\boldsymbol \Gamma= \mathrm{axl}\big(\boldsymbol Q_e^T\boldsymbol Q_{e,k}\big)\otimes \boldsymbol g^k\,,
\qquad\textrm{where} \qquad 
\boldsymbol \Gamma (\boldsymbol x):=  \boldsymbol \Gamma_\xi\big(\boldsymbol\Theta(\boldsymbol x)\big)
.
\end{equation}
In order to establish \eqref{e14,1}, we need to show that
\begin{equation}\label{e14,2}
\mathrm{axl}\big(\boldsymbol Q_e^T\boldsymbol Q_{e,k}\big)\otimes \boldsymbol g^k= \mathrm{axl}\big(\boldsymbol R_\xi^T\, \partial_{\xi_i}\boldsymbol R_\xi\big)\otimes \boldsymbol e_i\,.
\end{equation}
First, using the chain rule, we find
\begin{equation}\label{e14,3}
\boldsymbol Q_{e,k} = \dfrac{\partial}{\partial x_k}\,\, \boldsymbol R_\xi\big(\boldsymbol\Theta(\boldsymbol x)\big) = \dfrac{\partial\boldsymbol R_\xi }{\partial \xi_i}\,\, \dfrac{\partial \xi_i  }{\partial x_k}
 \qquad\Longrightarrow \qquad  
 \boldsymbol Q_e^T\boldsymbol Q_{e,k} =  \big(\boldsymbol R_\xi^T\, \partial_{\xi_i}\boldsymbol R_\xi\big)\, \dfrac{\partial \xi_i  }{\partial x_k}\;.
\end{equation}
Since $ \mathrm{axl}(\boldsymbol{A})= -\frac{1}{2}\,\boldsymbol{\epsilon}: \boldsymbol{A}$ for any skew-symmetric tensor $ \boldsymbol{A} $, where $ \boldsymbol{\epsilon} $ is the alternating Ricci third-order tensor (see e.g.\ \cite{Birsan-Neff-L57-2016,Birsan-Neff-L58-2017}), relation \eqref{e14,3} yields
\begin{equation}\label{e14,4}
\mathrm{axl}\big(\boldsymbol Q_e^T\boldsymbol Q_{e,k}\big) =  \mathrm{axl}\big(\boldsymbol R_\xi^T\, \partial_{\xi_i}\boldsymbol R_\xi\big)\, \dfrac{\partial \xi_i  }{\partial x_k}\;.
\end{equation}
Therefore, since $ \, \dfrac{\partial \xi_i  }{\partial x_k}\,=\boldsymbol e_i\cdot\boldsymbol g_k \; $, equation \eqref{e14,2} and thus \eqref{e14,1} holds.

Using the relation $ \,\boldsymbol Q_e= \overline{\boldsymbol R} \,\boldsymbol Q_0^T\, $ and the identity 
\[ \boldsymbol Q \,\mathrm{axl}(\boldsymbol A) = \mathrm{axl}\big(\boldsymbol Q\boldsymbol A\boldsymbol Q^T\big)\,,\] 
which is valid for any skew-symmetric tensor $ \boldsymbol{A} \in\mathfrak{so}(3)$ and any $\boldsymbol Q\in  \mathrm{SO}(3)$, we can write the relation \eqref{e14,1} in the alternative form 
\begin{equation}\label{e14,5}
    \boldsymbol \Gamma=  \boldsymbol Q_0\big[  \mathrm{axl}\big(\overline{\boldsymbol{R}}^T\overline{\boldsymbol{R}}_{,i}\big) -  \mathrm{axl}\big(\boldsymbol Q_0^T\boldsymbol Q_{0,i}\big) \big] \otimes \boldsymbol g^i\,.
\end{equation}
The wryness tensor defined in \eqref{e14,1} can also be written in the form
\begin{equation}\label{e14,6}
\boldsymbol \Gamma =  \boldsymbol \Gamma_e \big(\nabla_x \boldsymbol\Theta\big)^{-1} \,,
\qquad\textrm{where} \qquad 
\boldsymbol \Gamma_e:=  \mathrm{axl}\big(\boldsymbol Q_e^T\boldsymbol Q_{e,i}\big)\otimes  \boldsymbol e_i
\;.
\end{equation}
For a detailed discussion on various strain measures of non-linear micropolar continua we refer to the papers \cite{Pietraszkiewicz09,Birsan-Neff-L57-2016}.
\bigskip

We now turn to the constitutive relations.
We assume that the elastically stored energy density $W(\overline{\boldsymbol{E}}, \boldsymbol \Gamma)$ admits the additive split
\begin{equation}\label{e15}
    W(\overline{\boldsymbol{E}}, \boldsymbol \Gamma)=W_{\mathrm{mp}}(\overline{\boldsymbol{E}})+ W_{\mathrm{curv}}(  \boldsymbol \Gamma)
\end{equation}
into the elastic stretch (membrane) part $W_{\mathrm{mp}}$ and the curvature part $W_{\mathrm{curv}}$. For the elastic stretch part, we assume the form
\begin{equation}\label{e16}
    \begin{array}{rcl}
      W_{\mathrm{mp}}(\overline{\boldsymbol{E}} ) &=& \mu\,\|\,\mathrm{sym}\, \overline{\boldsymbol{E}}\,\|^2\, +  \,\mu_c \, \|\,\mathrm{skew}\, \overline{\boldsymbol{E}}\,\|^2\,   + \, \dfrac{\lambda}{2}\,\big(\mathrm{tr}\,\overline{\boldsymbol{E}}\,\big)^2\vspace{4pt}\\
        &=&  \mu\,\|\,\mathrm{dev_3\,sym}\, \overline{\boldsymbol{E}}\,\|^2\, +  \,\mu_c \, \|\,\mathrm{skew}\, \overline{\boldsymbol{E}}\,\|^2\,   + \, \dfrac{\kappa}{2}\,\big(\mathrm{tr}\,\overline{\boldsymbol{E}}\,\big)^2\,,
    \end{array}
\end{equation}
where $\lambda,\mu$ are the Lam\'e constants and $\kappa=\frac{3\lambda+2\mu}3\,$ is the bulk modulus of classical isotropic elasticity, while $\,\mu_c\ge 0$ is called the \emph{Cosserat couple modulus} \cite{Neff_zamm06,Neff_plate07_m3as};
here, we also employ the well-known operators
\begin{equation} 
\mathrm{sym}\,\boldsymbol{X}= \frac12 \,(\boldsymbol{X}+\boldsymbol{X}^T ) ,\qquad \mathrm{skew}\,\boldsymbol{X} = \frac12 \,(\boldsymbol{X}-\boldsymbol{X}^T ),\qquad
\,\mathrm{dev_3\,}\boldsymbol{X}=\boldsymbol{X}-\frac13(\mathrm{tr}\,\boldsymbol{X})\,\id_3\,,
\end{equation}
which represent the symmetric part, the skew-symmetric part and the deviatoric part, respectively, of any three-dimensional tensor $\boldsymbol{X}\,$. We also assume the standard restrictions
\begin{equation}\label{e17}
    \mu>0,\qquad \kappa>0\quad\mathrm{and}\quad  \mu_c\ge 0
\end{equation}
on the constitutive coefficients.

For the curvature part $ W_{\mathrm{curv}}(  \boldsymbol \Gamma)$ of the energy  density, we use the form \cite{Birsan-Neff-Ost_L56-2015}
\begin{equation}\label{e18}
    W_{\mathrm{curv}}( \boldsymbol{\Gamma}) = \mu\,L_c^2\,\Big(\,b_1\,\|\,\mathrm{dev_3\,sym}\, \boldsymbol{\Gamma}\|^2\, +  \,b_2\, \|\,\mathrm{skew}\, \boldsymbol{\Gamma}\|^2\,   + \, b_3\big(\mathrm{tr}\,\boldsymbol{\Gamma}\big)^2\,\Big)\,,
\end{equation}
where  $b_1\,,\, b_2 \,,\, b_3>0$ are dimensionless constitutive coefficients and the parameter $\,L_c>0\,$ introduces an internal length which is characteristic for the material and is responsible for size effects of the Cosserat model.
We notice that the energy density \eqref{e15} is a quadratic function, which corresponds to a physically linear response.

Assuming that no external body and surface forces and no external volume or surface couples are present, the deformation and the microrotation solve the geometrically nonlinear minimization problem
\begin{equation}\label{e19}
    I  =\dd\int_{\Omega_\xi}\widetilde W \big(\boldsymbol F_\xi ,{\boldsymbol{R}}_\xi\,\big)\, \mathrm dV\quad\to \quad \textrm{\ \ min\ \   w.r.t.\ \ } (\boldsymbol\varphi_\xi,{\boldsymbol{R}}_\xi\,)\, ,
\end{equation}
posed on $\Omega_\xi\,$, where $\widetilde W \big(\boldsymbol F_\xi ,{\boldsymbol{R}}_\xi\,\big)=  W(\overline{\boldsymbol{E}}, \boldsymbol \Gamma) $ is the energy density.
Making the change of variables $(\xi_1,\xi_2,\xi_3)=\boldsymbol \Theta(\boldsymbol x)  $ in the above integral we can write the total energy functional $I$ as the integral
\begin{equation}\label{e20}
     I  =\dd\int_{\Omega_h}  \big(\,W_{\mathrm{mp}}(\overline{\boldsymbol{E}})+  W_{\mathrm{curv}}(  \boldsymbol \Gamma)\,\big)\,\mathrm{det} \big[  \nabla_x  \boldsymbol\Theta(\boldsymbol x) \big]\, \mathrm d\boldsymbol x\,.
\end{equation}
over the Cartesian domain $\Omega_h\,$. Under the constitutive assumptions \eqref{e15}--\eqref{e18}, the existence of minimizers for the minimization problem \eqref{e19} for three-dimensional Cosserat models has been presented, e.g., in \cite{Neff_Edinb06,Birsan-Neff-Ost_L56-2015}, provided that, in addition, $\,\mu_c>0\,  $.

In the next sections, we will confine our attention to thin domains and deduce an improved Cosserat shell model by dimensional reduction.

\setcounter{equation}{0}

\section{The three--dimensional problem on a thin domain}

In the following, we assume that the parameter domain $\Omega_h\subset\R^3$ is a right cylinder of the form
$$\Omega_h=\left\{ (x_1,x_2,x_3) \,\Big|\,\, (x_1,x_2)\in\omega, \,\,\,-\dfrac{h}{2}\,< x_3<\, \dfrac{h}{2}\, \right\} =\,\,\dd\omega\,\times\left(-\frac{h}{2}\,\,,\,\frac{h}{2}\right),$$
where  $\omega\subset \mathbb{R}^2$ is a plane domain and the constant length $h>0$ is the thickness of the shell.
For shell--like bodies we assume  the  domain $\Omega_h $ to be \emph{thin}, i.e. the length $h$ to be \emph{small} in comparison with its diameter.

Furthermore, for our purpose we assume that the parametric representation $\boldsymbol\Theta$ has the form
\begin{equation}\label{e21}
    \boldsymbol\Theta(\boldsymbol x )=\boldsymbol y_0(x_1,x_2)+x_3\, \boldsymbol n_0(x_1,x_2), \qquad\boldsymbol n_0=\dfrac{\boldsymbol y_{0,1}\times \boldsymbol y_{0,2}}{\|\boldsymbol y_{0,1}\times \boldsymbol y_{0,2}\|}\,,
\end{equation}
i.e.\ that $\boldsymbol\Theta$ maps the midsurface $\omega$ of the  parameter domain $\Omega_h$ onto the midsurface  $\omega_\xi:=\boldsymbol y_0(\omega)$ of the reference configuration $\Omega_\xi\,$, while $\boldsymbol n_0$ is the unit normal vector to $\omega_\xi\,$. The special form \eqref{e21} is a classical representation in shell theory, see e.g., \cite{Ciarlet00,Libai98,Pietraszkiewicz-book04}.

\subsection{Prerequisites from classical differential geometry of surfaces}\label{sectGeo}

In preparation for the dimensional reduction, we need to state a number of well-known formulas from the differential geometry of surfaces in $\R^3$ (applied to the reference midsurface $\omega_\xi\,$) as well as additional notational conventions. The midsurface $\omega_\xi\,$ of the reference configuration $\Omega_\xi\,$ admits the parametric representation $\boldsymbol y_0(x_1,x_2)$, with $(x_1,x_2)\in\omega$. We introduce the covariant base vectors $\boldsymbol a_1\,,\boldsymbol a_2\,$ and the contravariant base vectors $\boldsymbol a^1\,,\boldsymbol a^2\,$ in the tangent plane by
\begin{equation}\label{e22}
    \boldsymbol a_\alpha:=\,\dfrac{\partial \boldsymbol y_0}{\partial x_\alpha}\,=\boldsymbol y_{0,\alpha}\,\,,\qquad \boldsymbol a^\beta\cdot\boldsymbol a_\alpha=\delta_\alpha^\beta\,\,,\qquad \alpha,\beta=1,2\,,
\end{equation}
and let $\boldsymbol a_3\,=\boldsymbol a^3=\boldsymbol n_0\,$.
The \emph{first fundamental tensor} $\boldsymbol a$ of the surface $\omega_\xi$ is
\begin{equation}\label{e23}
\boldsymbol{a}:=\boldsymbol{a}_\alpha\otimes \boldsymbol{a}^\alpha=
a_{\alpha\beta}\boldsymbol{a}^\alpha\otimes \boldsymbol{a}^\beta= a^{\alpha\beta}\boldsymbol{a}_\alpha\otimes \boldsymbol{a}_\beta \quad\mathrm{with}  \quad
a_{\alpha\beta}=\boldsymbol{a}_\alpha\cdot \boldsymbol{a}_\beta\,,\quad a^{\alpha\beta}=\boldsymbol{a}^\alpha\cdot \boldsymbol{a}^\beta .
\end{equation}
We denote by
\begin{equation*}
a(x_1,x_2):=\sqrt{\mathrm{det}\, \big( a_{\alpha\beta} \big){}_{2\times 2}}\,>0
\end{equation*}
the quotient between the area element of $ \omega_\xi $ and the area element of $ \omega $.
Let the operator $\mathrm{Grad}_s\,$ denote the \emph{gradient on the surface}, given for any field $\,\boldsymbol f\,$ on $\omega_\xi\,$ by
$$\mathrm{Grad}_s\,\boldsymbol f\,:=\, \dfrac{\partial\boldsymbol f}{\partial x_\alpha}\,\otimes \boldsymbol a^\alpha= \boldsymbol f,_{\alpha}\otimes \boldsymbol a^\alpha.$$
Then $ \;\boldsymbol{a}= \text{Grad}_s\,\boldsymbol{y}_0 $ . The \emph{second fundamental tensor} $\boldsymbol b$ of the surface is expressed by
\begin{equation}\label{e24}
\begin{array}{c}
  \boldsymbol{b}:= -\text{Grad}_s\,\boldsymbol{n}_0=-  \boldsymbol{n}_{0,\alpha}\otimes\boldsymbol{a}^\alpha= b_{\alpha\beta}\,\boldsymbol{a}^\alpha\otimes \boldsymbol{a}^\beta=b^\alpha_\beta\,\boldsymbol{a}_\alpha\otimes \boldsymbol{a}^\beta,
  \vspace{4pt}\\
    b_{\alpha\beta}=-  \boldsymbol{n}_{0,\beta}\cdot\boldsymbol{a}_\alpha= b_{\beta\alpha}\,,\qquad b^{\alpha}_{\beta}=-  \boldsymbol{n}_{0,\beta}\cdot \boldsymbol{a}^\alpha\,;
\end{array}
\end{equation}
recall that the tensor $\boldsymbol b$ is symmetric. We also employ the usual notations for the \emph{mean curvature} $H$ and the \emph{Gau\ss{} curvature} $K$ of the surface
\begin{equation}\label{e25}
    2H:=\mathrm{tr}\,\boldsymbol b=b^{\alpha}_{\alpha}\,,\qquad K:=\mathrm{det}\,\boldsymbol b=\mathrm{det}\, \big( b^{\alpha}_{\beta} \big){}_{2\times 2}= \dfrac{1}{2}\,\big[ \big(\mathrm{tr}\,\boldsymbol b \,\big)^2 - \mathrm{tr}\big(\boldsymbol b^2\big)\big].
\end{equation}
Then the Cayley-Hamilton theorem, applied to the tensor $\boldsymbol b\,$, yields
\begin{equation}\label{e26}
    \boldsymbol b^2-2H\cdot \boldsymbol b+K\cdot\boldsymbol a=\boldsymbol 0.
\end{equation}
Finally, the so-called \emph{alternator tensor} $\boldsymbol c$ of the surface \cite{Zhilin06} is given by
\begin{equation}\label{e27}
    \boldsymbol c:=-\boldsymbol n_0\times \boldsymbol a=-\boldsymbol a\times\boldsymbol n_0=  \dfrac{1}{a}\,\,\epsilon_{\alpha\beta}\,\boldsymbol{a}_\alpha\otimes \boldsymbol{a}_\beta = a\,\,\epsilon_{\alpha\beta}\,\boldsymbol{a}^\alpha\otimes \boldsymbol{a}^\beta,
\end{equation}
where $\epsilon_{\alpha\beta}\,$ is the two-dimensional alternator with $\epsilon_{12}=-\epsilon_{21}=1\,,\,\epsilon_{11}=\epsilon_{22}=0$. Note that the tensor $\boldsymbol c$ is antisymmetric and satisfies $\, \boldsymbol c^2=-\boldsymbol a$. Furthermore, the tensors $\boldsymbol a\,$, $\boldsymbol b\,$, and $\boldsymbol c$ defined above are \emph{planar}, i.e.\ tensors in the tangent plane of the surface, with $\boldsymbol a$ being the identity tensor in the tangent plane. The notation introduced above will be used throughout the rest of the article.

\subsection{Useful relations for the  gradient of the mapping $\,\boldsymbol \Theta$}\label{sectDefGrad}

Next, we will write the expressions of the base vectors $\boldsymbol g_i\,$, $\boldsymbol g^i$, the gradient $\nabla_x \boldsymbol\Theta$ and the inverse $\big[\nabla_x \boldsymbol\Theta\big]^{-1}$ corresponding to the special form of the mapping $\boldsymbol \Theta$ given by \eqref{e21}. From \eqref{e21}--\eqref{e26} we get the known relations (see e.g., \cite{Pietraszkiewicz14})
\begin{equation}\label{e28}
    \begin{array}{l}
      \boldsymbol g_\alpha= \boldsymbol\Theta_{,\alpha}= \boldsymbol a_\alpha+ x_3\,\boldsymbol n_{0,\alpha}=\big(\delta_\alpha^\beta- x_3\,b_\alpha^\beta\,\big)\,\boldsymbol a_\beta=\big(\boldsymbol a-x_3\,\boldsymbol b\big)\,\boldsymbol a_\alpha\,,\vspace{3pt}\\
      \boldsymbol g_3= \boldsymbol g^3= \boldsymbol n_0\,, \vspace{3pt}\\
      \boldsymbol g^\alpha=\,\dfrac{1}{b(x_3)}\,\big[\, \delta^\alpha_\beta+ x_3\big(b^\alpha_\beta - 2H\, \delta^\alpha_\beta\big)  \big]\boldsymbol a^\beta =\,\dfrac{1}{b(x_3)}\,\big[ \boldsymbol a+ x_3\big(\boldsymbol b - 2H\, \boldsymbol a\big)  \big]\boldsymbol a^\alpha
    \end{array}
\end{equation}
with $
b(x_3):= 1-2H\,x_3+K\,x_3^2\,.$

For the gradient of   $\boldsymbol \Theta$ we obtain from \eqref{e5}$_{1,2}$ and \eqref{e28}$_{1,2}$
\begin{equation}\label{e29}
    \nabla_x \boldsymbol\Theta\,=\, \boldsymbol g_i\otimes\boldsymbol e_i= \big(\delta^\alpha_\beta- x_3\,b^\alpha_\beta\,\big)\,\boldsymbol a_\alpha\otimes \boldsymbol e_\beta + \boldsymbol n_0\otimes \boldsymbol e_3
\end{equation}
and
\begin{equation}\label{e30}
    \mathrm{det}\,\big( \nabla_x \boldsymbol\Theta \big)= a(x_1,x_2)\cdot\big[ 1-2H\,x_3+K\,x_3^2 \,\big]= a(x_1,x_2)\cdot b(x_3)\,>0.
\end{equation}
We observe that $b(x_3)\,>0$ for $x_3\in(-\frac{h}{2}\,,\,\frac{h}{2})$ and sufficiently small $h$. Indeed, we assume that in general 
$$h\ll R_1 \qquad\textrm{and}\qquad h\ll R_2\,,$$ 
where $R_\alpha= |K_\alpha|^{-1}$, $ K_\alpha $ are the principal curvatures of the surface $\omega_\xi\,$, with $2H=K_1+K_2\,,\,K=K_1K_2\,$, and thus
$$b(x_3)=(1-x_3K_1)(1-x_3K_2)>0.$$
Similarly, from \eqref{e5}$_{2}$ and \eqref{e28} we obtain
\begin{equation}\label{e31}
\begin{array}{cl}
  \big[\nabla_x \boldsymbol\Theta\big]^{-1} & = \boldsymbol e_i\otimes\boldsymbol g^i= \,\dfrac{1}{b(x_3)}\,\big[\, \delta^\alpha_\beta+ x_3\big(b^\alpha_\beta - 2H\, \delta^\alpha_\beta\big)  \big]\boldsymbol e_\alpha\otimes\boldsymbol a^\beta +\boldsymbol e_3\otimes \boldsymbol n_0 \\
    & =  \,\dfrac{1}{b(x_3)}\,\Big\{ \big( \boldsymbol e_i\otimes\boldsymbol a^i \big) +x_3   \big[ \big(b^\alpha_\beta - 2H \delta^\alpha_\beta\big) \boldsymbol e_\alpha\otimes\boldsymbol a^\beta - 2H\,\boldsymbol e_3\otimes \boldsymbol n_0\big]+ x_3^2\,K\boldsymbol e_3\otimes \boldsymbol n_0\Big\}.
\end{array}
\end{equation}
By virtue of \eqref{e30} and \eqref{e31},
\begin{align}\label{e32}
      \mathrm{Cof}\big( \nabla_x \boldsymbol \Theta \big)  &= \mathrm{det}\big( \nabla_x \boldsymbol \Theta \big)\cdot  \big( \nabla_x \boldsymbol \Theta \big)^{-T}  \vspace{3pt}\\
      & =a(x_1,x_2)\Big\{
      \big( \boldsymbol a^i\otimes\boldsymbol e_i \big) +x_3   \big[ \big(b^\alpha_\beta - 2H \delta^\alpha_\beta\big)\boldsymbol a^\beta \otimes \boldsymbol e_\alpha - 2H\, \boldsymbol n_0\otimes \boldsymbol e_3\big]+ x_3^2\,K\boldsymbol n_0\otimes \boldsymbol e_3  \Big\}.\notag
\end{align}
From \eqref{e29}, we also deduce
\begin{equation}\label{e33}
    \big(\nabla_x \boldsymbol\Theta\big)^T\big(\nabla_x \boldsymbol\Theta\big) = \big(\delta_{\alpha\gamma}- x_3\,b_{\alpha\gamma}\,\big) \big(\delta^\gamma_\beta- x_3\,b^\gamma_\beta\,\big)\,\boldsymbol e_\alpha\otimes \boldsymbol e_\beta + \boldsymbol e_3\otimes \boldsymbol e_3\,.
\end{equation}

Recall from \eqref{e7} that $\boldsymbol d_3^0=  \boldsymbol Q_0\boldsymbol e_3\,$ for $\boldsymbol Q_0 = \polar(\nabla_x \boldsymbol \Theta) $. We will show that
\begin{equation}\label{e34}
    \boldsymbol d_3^0=\boldsymbol n_0\,.
\end{equation}
To this aim, we first need to show that the symmetric positive definite tensor $\boldsymbol U_0$, as defined by \eqref{e6}, has the form
\begin{equation}\label{e35}
    \boldsymbol U_0(\boldsymbol x)=u_{\alpha\beta}(\boldsymbol x)\boldsymbol e_\alpha\otimes \boldsymbol e_\beta+ \boldsymbol e_3\otimes \boldsymbol e_3\,,
\end{equation}
with $u_{\alpha\beta}=u_{ \beta\alpha}$ and $\mathrm{det}(u_{\alpha\beta})_{2\times 2}>0$. Indeed, due to \eqref{e6} and \eqref{e33}, we find
\begin{equation}\label{e36}
    \boldsymbol U_0^2= \big(\nabla_x \boldsymbol\Theta\big)^T\big(\nabla_x \boldsymbol\Theta\big) = v_{\alpha\beta} \,\boldsymbol e_\alpha\otimes \boldsymbol e_\beta+ \boldsymbol e_3\otimes \boldsymbol e_3\,,
\end{equation}
where $v_{\alpha\beta} =  \delta_{\alpha\beta}- 2x_3\,b_{\alpha\beta}+ x_3^2\,b_{\alpha\gamma}\,b^\gamma_\beta\, $. Let $\lambda_1,\lambda_2,\lambda_3>0$ be the eigenvalues of the tensor $\boldsymbol U_0\,$. From \eqref{e36} we see that $\boldsymbol U_0^2$ admits the eigenvalue 1 (with corresponding eigenvector $\boldsymbol e_3$). Consequently, $\lambda_3=1$ is an eigenvalue of $\boldsymbol U_0$ and thus
$\mathrm{tr}\,\boldsymbol U_0= \lambda_1+\lambda_2+1$, $\mathrm{tr\,(Cof\,}\boldsymbol U_0)= \lambda_1\lambda_2+\lambda_1+\lambda_2\,$, and $\mathrm{det\,}\boldsymbol U_0=  \lambda_1\lambda_2\,$. Hence, the characteristic equation of $\boldsymbol U_0$ is
\begin{align}
     & \boldsymbol U_0^3-(\lambda_1+\lambda_2+1)\,\boldsymbol U_0^2+(\lambda_1\lambda_2+\lambda_1+\lambda_2)\,\boldsymbol U_0-(\lambda_1\lambda_2)\,\id_3= \boldsymbol 0\qquad\  \Leftrightarrow \vspace{4pt}\notag\\
     & \boldsymbol U_0\,\big[\boldsymbol U_0^2 +(\lambda_1\lambda_2+\lambda_1+\lambda_2)\,\id_3\big] = (\lambda_1+\lambda_2+1)\,\boldsymbol U_0^2+(\lambda_1\lambda_2)\,\id_3 \qquad \Leftrightarrow 
 \notag\\&\label{e36,5}
      \boldsymbol U_0\,\big[\big(  v_{\alpha\beta} +(\lambda_1\lambda_2+\lambda_1+\lambda_2) \delta_{\alpha\beta}\big)\,\boldsymbol e_\alpha\otimes \boldsymbol e_\beta+ (\lambda_1\lambda_2+\lambda_1+\lambda_2+1)  \,\boldsymbol e_3\otimes \boldsymbol e_3\big]
       \vspace{4pt}\\
       &\qquad\qquad = \big[(\lambda_1+\lambda_2+1) v_{\alpha\beta} +(\lambda_1\lambda_2) \delta_{\alpha\beta}\big]\,\boldsymbol e_\alpha\otimes \boldsymbol e_\beta+ (\lambda_1\lambda_2+\lambda_1+\lambda_2+1)  \,\boldsymbol e_3\otimes \boldsymbol e_3\,.\notag
    \end{align}
Since the eigenvalues of the matrix $(v_{\alpha\beta})_{2\times 2}$ are $\lambda_1^2\,$ and $\,\lambda_2^2\,$, we find
\begin{equation*}
    \begin{array}{cl}
      \mathrm{det\,}\big( v_{\alpha\beta} +(\lambda_1\lambda_2+\lambda_1+\lambda_2) \delta_{\alpha\beta}\big)_{2\times 2} & =(\lambda_1^2+\lambda_1\lambda_2+\lambda_1+\lambda_2) (\lambda_2^2+\lambda_1\lambda_2+\lambda_1+\lambda_2)  \\
       & =(\lambda_1+\lambda_2)^2(\lambda_1 +1)( \lambda_2+1)>0.
    \end{array}
\end{equation*}
Hence, the matrix $\big( v_{\alpha\beta} +(\lambda_1\lambda_2+\lambda_1+\lambda_2) \delta_{\alpha\beta}\big)_{2\times 2}$ is invertible and we denote its inverse by $\big( w_{\alpha\beta}\big)_{2\times 2}\,$. Then \eqref{e36,5} implies
$$ \boldsymbol U_0=  \Big\{ \big[(\lambda_1+\lambda_2+1) v_{\alpha\beta} +(\lambda_1\lambda_2) \delta_{\alpha\beta}\big]\,\boldsymbol e_\alpha\otimes \boldsymbol e_\beta+ (\lambda_1 +1)( \lambda_2+1)  \,\boldsymbol e_3\otimes \boldsymbol e_3\Big\}\,\Big[   w_{\alpha\beta}  \,\boldsymbol e_\alpha\otimes \boldsymbol e_\beta + \,\dfrac{\boldsymbol e_3\otimes \boldsymbol e_3}{(\lambda_1+1)(\lambda_2+1)}  \,\Big]\,,$$
and by multiplication on the right-hand side we see that $\boldsymbol U_0\,$ does indeed have the form \eqref{e35}, which in turn implies
$$\boldsymbol U_0\,\boldsymbol e_3=\boldsymbol e_3\qquad\Rightarrow\qquad \boldsymbol e_3 = \boldsymbol U_0^{-1}\,\boldsymbol e_3\,.$$
Therefore, and in view of \eqref{e6}, \eqref{e7} and \eqref{e29},
$$\boldsymbol d_3^0=\boldsymbol Q_0\,\boldsymbol e_3=   \big(\nabla_x \boldsymbol\Theta\big)\boldsymbol U_0^{-1}\, \boldsymbol e_3 =  \big(\nabla_x \boldsymbol\Theta\big)\boldsymbol e_3  = \boldsymbol n_0\,,$$
which establishes the relation \eqref{e34}. This means that the initial director $\boldsymbol d_3^0$ is chosen along the normal to the reference midsurface (the ``material filament'' of the shell), while $\{\boldsymbol d_1^0\,,\,\boldsymbol d_2^0\,\}$ is an orthonormal basis in the tangent plane of $\omega_\xi\,$. Then, we can express the tensors $\boldsymbol a$ and $\boldsymbol c$  defined by \eqref{e23} and \eqref{e27}   in the alternative forms
\begin{equation}\label{e37}
    \boldsymbol{a}=
    \boldsymbol{d}^0_\alpha\otimes \boldsymbol{d}^0_\alpha\,,\qquad  \boldsymbol{c}=  \epsilon_{\alpha\beta}\,\boldsymbol{d}^0_\alpha\otimes \boldsymbol{d}^0_\beta\,.
\end{equation}

In  the current configuration $\Omega_c$ the director $\boldsymbol d_3$ is no longer orthogonal to the deformed surface $\omega_c:=\boldsymbol \varphi(\omega)$ and the vectors $\boldsymbol d_1\,,\, \boldsymbol d_2\,$ are not tangent to this surface. The deviation of the director $\boldsymbol d_3$ from the normal vector to $\omega_c$ describes the \emph{transverse shear deformation} of shells. Moreover, the rotations of $\{\boldsymbol d_1\,,\,\boldsymbol d_2\,\}$ about the  director $\boldsymbol d_3$ describe the so-called \emph{drilling rotations} in shells (see \cite{Birsan-Neff-L54-2014}).

\subsection{Stress tensors of Piola--Kirchhoff type}\label{sectPiolaKirch}

We consider the analogue of the second Piola--Kirchhoff stress tensor  from   classical elasticity theory, given by the derivative
\begin{equation}\label{e38}
    \boldsymbol{S}_2\big( \overline{\boldsymbol E} \big):=
    D_{\overline{\boldsymbol E}}\,\,W_{\rm mp}\big( \overline{\boldsymbol E} \big)
\end{equation}
and the analogue of the first Piola--Kirchhoff stress tensor  given by
\begin{equation}\label{e39}
    \boldsymbol{S}_1 \big(\boldsymbol F ,\overline{\boldsymbol{R}}\,\big):=
    D_{ \boldsymbol F}\,\widetilde W_{\mathrm{mp}} \big(\boldsymbol F ,\overline{\boldsymbol{R}}\,\big),
\end{equation}
where $\widetilde W_{\mathrm{mp}} \big(\boldsymbol F ,\overline{\boldsymbol{R}}\,\big)= W_{\rm mp}\big( \overline{\boldsymbol E} \big)$ is the elastic stretch part of the energy density, expressed as a function of the deformation gradient $\boldsymbol F$ and the total microrotation $\overline{\boldsymbol{R}}$. Note that
\begin{equation}\label{e40}
    D_{ \boldsymbol F}\,\widetilde W_{\mathrm{mp}} \big(\boldsymbol F ,\overline{\boldsymbol{R}}\,\big)=
\boldsymbol Q_e\,
\big[D_{\overline{\boldsymbol E}}\,\,W_{\rm mp}\big( \overline{\boldsymbol E} \big)\big]\big(\nabla_x \boldsymbol\Theta\big)^{-T}
\quad\Rightarrow\quad
 \boldsymbol{S}_1 \big(\boldsymbol F ,\overline{\boldsymbol{R}}\,\big) =
\boldsymbol Q_e\,
\boldsymbol{S}_2\big( \overline{\boldsymbol E} \big)\big(\nabla_x \boldsymbol\Theta\big)^{-T}.
\end{equation}
Then, in view of \eqref{e16} and \eqref{e38} we find
\begin{equation}\label{e41}
    \begin{array}{rcl}
      \boldsymbol{S}_2\big( \overline{\boldsymbol E} \big) & = &
       2\mu\, \mathrm{sym}\, \overline{\boldsymbol{E}}\,  +  \,2\mu_c \,  \mathrm{skew}\, \overline{\boldsymbol{E}}\, \,   + \, \lambda \,\big( \mathrm{tr}\,\overline{\boldsymbol{E}}\big)\,\id_3\,\,,
      \vspace{4pt}\\
     \boldsymbol{S}_1 \big(\boldsymbol F ,\overline{\boldsymbol{R}}\,\big) & =  &
        \boldsymbol Q_e\,\Big[   2\mu\, \mathrm{sym}\, \overline{\boldsymbol{E}}\,  +  \,2\mu_c \,  \mathrm{skew}\, \overline{\boldsymbol{E}}\, \,   + \, \lambda \, \big(\mathrm{tr}\,\overline{\boldsymbol{E}}\big)\,\id_3\, \Big] \big(\nabla_x \boldsymbol\Theta\big)^{-T}.
    \end{array}
\end{equation}
We note that the tensor $ \boldsymbol{S}_2 $ is not symmetric in general, with its skew-symmetric part being governed by the Cosserat couple modulus $ \mu_c\, $.
 
As usual in the theory of shells, we shall assume that the stress vectors on the upper and lower surfaces (major surfaces) of the shell have null normal components, i.e.
\begin{equation}\label{e42}
    \big[ \boldsymbol{S}_1 \big(\boldsymbol F ,\overline{\boldsymbol{R}}\,\big)\cdot \boldsymbol e_3 \big]\cdot\boldsymbol d_3=0\qquad\mathrm{for}\quad x_3=\pm\,\dfrac{h}{2}\,\,.
\end{equation}
In view of \eqref{e28} we see indeed that  $\,\,\boldsymbol g_\alpha\cdot\boldsymbol n_0=0\,$, which means that $\boldsymbol n_0$ is also normal to the major surfaces of the shell (i.e. the upper and lower surfaces, characterized by $x_3=\pm\,\frac{h}{2}\,$). Then the outward unit normals to the shell boundary are $\boldsymbol n_0$ for $x_3= \,\frac{h}{2}\,$ and $-\boldsymbol n_0$ for $x_3=-\,\frac{h}{2}\,$, respectively.

A similar condition to \eqref{e42} was also employed in the derivation of the Koiter shell model from three-dimensional nonlinear elasticity (see e.g., \cite{Steigmann12,Steigmann13}).

Using \eqref{e40} and $\,\boldsymbol Q_e^T\boldsymbol d_3=\boldsymbol d_3^0=\boldsymbol n_0\,\,$ we can write \eqref{e42} in the form
\begin{align}
    \big[ \boldsymbol{S}_2\big( \overline{\boldsymbol E} \big)\big(\nabla_x \boldsymbol\Theta\big)^{-T}\cdot \boldsymbol e_3 \big]\cdot\boldsymbol n_0&=0\qquad\mathrm{for}\quad x_3=\pm\,\dfrac{h}{2}\,\,, \quad\mbox{i.e.}\notag
\\
\big[ \boldsymbol{S}_2\big( \overline{\boldsymbol E} \big)\cdot \boldsymbol n_0 \big]\cdot\boldsymbol n_0&=0\qquad\mathrm{for}\quad x_3=\pm\,\dfrac{h}{2}\,\,,\label{e43}
\end{align}
since $ \big(\nabla_x \boldsymbol\Theta\big)^{-T}\cdot \boldsymbol e_3 = \boldsymbol n_0 $ from \eqref{e31}.
A simplified approximate form of \eqref{e43} can be obtained in the limit as $h\to 0$. Indeed, if we denote by $f$ the function with
\begin{equation}\label{e44}
    f(z):= \Big\{\big[ \boldsymbol{S}_2\big( \overline{\boldsymbol E} \big)\,\big]_{x_3=z} \cdot\boldsymbol n_0 \Big\}\cdot\boldsymbol n_0 \qquad\mathrm{for}\quad z\in\Big[\,-\dfrac{h}{2}\,\,,\dfrac{h}{2}\,\Big],
\end{equation}
then from the Taylor expansion of $f(z)$ about $z=0$ we find
\begin{equation}\label{e45}
    \begin{array}{l}
      f\Big(  \dfrac{h}{2}\Big)+ f\Big(   \dfrac{-h}{2}\Big) = 2f(0)+O(h^2), \qquad
      f\Big(  \dfrac{h}{2}\Big)- f\Big( \dfrac{-h}{2}\Big) = hf'(0)+O(h^3),
\end{array}
\end{equation}
where
\begin{equation}\label{e45,5}     
f(0)= \Big\{\big[ \boldsymbol{S}_2\big( \overline{\boldsymbol E} \big)\,\big]_{x_3=0} \cdot\boldsymbol n_0 \Big\}\cdot\boldsymbol n_0\,,\qquad
       f'(0)= \Big\{\Big[ \,\dfrac{\partial}{\partial x_3}\,\boldsymbol{S}_2\big( \overline{\boldsymbol E} \big)\,\Big]_{x_3=0} \cdot\boldsymbol n_0 \Big\}\cdot\boldsymbol n_0\,.
\end{equation}
In view of \eqref{e43} and \eqref{e44}, we find $  f\big(  \frac{h}{2}\big)= f\big(   \frac{-h}{2}\big) =0$, and  the relations \eqref{e45} yield $f(0)=O(h^2)$ and $f'(0)=O(h^2) $. In the limit as $h\to 0$, the conditions \eqref{e43} are therefore approximated by $f(0)=0$ and $f'(0)=0 $\,, i.e.\ (in view of \eqref{e45,5})
\begin{equation}\label{e46}
    \Big\{\big[ \boldsymbol{S}_2\big( \overline{\boldsymbol E} \big)\,\big]_{x_3=0} \cdot\boldsymbol n_0 \Big\}\cdot\boldsymbol n_0=0\,,\qquad
        \Big\{\Big[ \,\dfrac{\partial}{\partial x_3}\,\boldsymbol{S}_2\big( \overline{\boldsymbol E} \big)\,\Big]_{x_3=0} \cdot\boldsymbol n_0 \Big\}\cdot\boldsymbol n_0=0\,.
\end{equation}
These relations will be used for the dimensional reduction procedure in the next section.

\setcounter{equation}{0}

\section{The 8--parameter ansatz}\label{Sect8Param}

In the following, we want to find a reasonable approximation $(\boldsymbol \varphi_s,\overline{\boldsymbol R}_{s})$ of the functions  $(\boldsymbol \varphi,\overline{\boldsymbol R})$ involving only two-dimensional quantities and show that this approximation is appropriate for shell-like bodies (here the subscript $s$ stands for shell).

We assume firstly that the total microrotation $\overline{\boldsymbol R}_s:\Omega_h\rightarrow \text{\rm{SO}}(3)$
in the thin shell do not depend on the thickness variable $x_3\,$, i.e. we set
\begin{equation}\label{e47}
        \overline{\boldsymbol R}_s(x_1,x_2)  = \overline{\boldsymbol R} (\boldsymbol x),
\end{equation}
which is in line with the assumed thinness and material homogeneity of the structure.
In the reference configuration, we similarly consider
\begin{equation}\label{e48}
    \boldsymbol Q_0 (\boldsymbol x)  =    \boldsymbol Q_0(x_1,x_2).
\end{equation}
In other words, the directors $\,\boldsymbol d_i^0$ and $\boldsymbol d_j$ are assumed to dependent only on the midsurface coordinates $(x_1,x_2)$:
\begin{equation}\label{e49}
 \boldsymbol d_i^0 (\boldsymbol x)  =    \boldsymbol d_i^0(x_1,x_2),\qquad
    \boldsymbol d_j (\boldsymbol x)  =    \boldsymbol d_j(x_1,x_2),\qquad i,j=1,2,3.
\end{equation}
For the elastic microrotation we obtain from \eqref{e47} and \eqref{e48}:
\begin{equation}\label{e50}
     \boldsymbol Q_e  =    \overline{\boldsymbol R}_s(x_1,x_2)\,\boldsymbol Q_0^T(x_1,x_2) = \boldsymbol Q_e(x_1,x_2).
\end{equation}

In the engineering shell community it is well known \cite{Chernykh80,Schmidt85,Pietraszkiewicz85} that
the ansatz for the deformation over the thickness should be at least  quadratic
in order to avoid the so called
{\it Poisson thickness locking}
and to fully capture the three-dimensional kinematics without artificial
modification of the material laws; see the detailed discussion of this point in \cite{Ramm00}
and compare with \cite{Ramm92,Ramm94,Ramm96,Ramm97,Sansour98c}.
We consider therefore the following 8-parameter quadratic ansatz in the thickness direction for the reconstructed total deformation $\boldsymbol \varphi_s:\Omega_h\subset \mathbb{R}^3\rightarrow \mathbb{R}^3$ of the shell-like structure
\begin{equation}\label{e51}
    \boldsymbol \varphi_s(\boldsymbol x)=\boldsymbol m(x_1,x_2)+\bigg(x_3\,\varrho_m(x_1,x_2)+\dd\frac{x_3^2}{2}\,\varrho_b(x_1,x_2)\bigg) \boldsymbol d_3(x_1,x_2) ,
\end{equation}
where $\boldsymbol d_3= \overline{\boldsymbol R}_s(x_1,x_2)\boldsymbol e_3= {\boldsymbol Q}_e(x_1,x_2)\boldsymbol n_0\,$ is the third director and  $\,\boldsymbol m: \omega\subset\mathbb{R}^2\rightarrow\mathbb{R}^3$ takes on the role of the deformation of the midsurface of
the shell viewed as a parametrized surface.   The yet indeterminate functions $\varrho_m,\,\varrho_b:\omega\subset \mathbb{R}^2\rightarrow \mathbb{R}$ allow in principal for symmetric thickness stretch  ($\varrho_m\neq1$) and asymmetric thickness stretch ($\varrho_b\neq 0$) about the midsurface.

We can now explicitly express the deformation gradient and strain measures corresponding to the assumed form of the deformation field \eqref{e51} and microrotation \eqref{e47}. In view of the relations \eqref{e11,25} and   \eqref{e28}, the (reconstructed) deformation gradient has the form
\begin{equation}\label{e52}
\begin{array}{crl}
 \boldsymbol F_s & := & \nabla_\xi\, \boldsymbol \varphi_s= \dfrac{\partial\boldsymbol\varphi_s }{\partial x_i}\,\,\otimes \boldsymbol g^i=
    \boldsymbol \varphi_{s,\alpha}\otimes\, \boldsymbol g^\alpha +
    \boldsymbol \varphi_{s,3}\otimes\, \boldsymbol n_0 \\
  & = & \dfrac{1}{b(x_3)}\,\Big[ \boldsymbol m_{,\alpha}+ x_3\big(\varrho_m \boldsymbol d_3 \big)_{,\alpha}+  \dfrac{x_3^2}{2}\,\big( \varrho_b \boldsymbol d_3\big)_{,\alpha} \Big] \otimes \boldsymbol a^\alpha\big[ \boldsymbol a+x_3(\boldsymbol b-2H\boldsymbol a) \big] +(\varrho_m+x_3\varrho_b)\boldsymbol d_3\otimes\boldsymbol n_0
  \\
  & = & \dfrac{1}{b(x_3)}\,\Big[ \mathrm{Grad}_s\boldsymbol m + x_3\mathrm{Grad}_s\big(\varrho_m \boldsymbol d_3 \big) +  \dfrac{x_3^2}{2}\,\mathrm{Grad}_s\big( \varrho_b \boldsymbol d_3\big)  \Big]
  \big[ \boldsymbol a+x_3(\boldsymbol b-2H\boldsymbol a) \big] \vspace{4pt}\\
  &  &+(\varrho_m+x_3\varrho_b)\boldsymbol d_3\otimes\boldsymbol n_0\,.
\end{array}
\end{equation}
Then from \eqref{e13}, we find
\begin{equation}\label{e53}
\begin{array}{crl}
  \boldsymbol{E}_s & := &\boldsymbol Q_e^T \boldsymbol F_s-\id_3\\
  &  = &
   \dfrac{1}{b(x_3)}\,\Big[ \boldsymbol Q_e^T\mathrm{Grad}_s\boldsymbol m + x_3\,\boldsymbol Q_e^T\mathrm{Grad}_s\big(\varrho_m \boldsymbol d_3 \big) +  \dfrac{x_3^2}{2}\,\boldsymbol Q_e^T\,\mathrm{Grad}_s\big( \varrho_b \boldsymbol d_3\big)  \Big]
  \big[ \boldsymbol a+x_3(\boldsymbol b-2H\boldsymbol a) \big]   \vspace{4pt}\\
  & & +(\varrho_m+x_3\varrho_b)\boldsymbol n_0\otimes\boldsymbol n_0 -\id_3
\end{array}
\end{equation}
and from \eqref{e14,1} and \eqref{e50} we obtain
\begin{equation}\label{e54}
\begin{array}{rl}
  \boldsymbol \Gamma_s:= & \mathrm{axl}\big(\boldsymbol Q_e^T\boldsymbol Q_{e,i}\big)\otimes \boldsymbol g^i=
     \mathrm{axl}\big(\boldsymbol Q_e^T\boldsymbol Q_{e,\alpha}\big)\otimes \boldsymbol g^\alpha \vspace{4pt}\\
   = & \dfrac{1}{b(x_3)}\,\, \mathrm{axl}\big(\boldsymbol Q_e^T\boldsymbol Q_{e,\alpha}\big)\otimes
   \boldsymbol a^\alpha\big[ \boldsymbol a+ x_3\big(\boldsymbol b - 2H\, \boldsymbol a\big)  \big].
\end{array}
\end{equation}

Next, we want to express the above tensors \eqref{e53} and \eqref{e54} with the help of  strain measures  used in the general nonlinear shell theory \cite{Eremeyev06}.  Therefore, we introduce the \emph{elastic shell strain tensor} $\boldsymbol E^e$ and the \emph{elastic shell bending--curvature tensor} $\boldsymbol K^e$, which are tensor fields on the surface $\omega_\xi\,$ defined by \cite{Libai98,Pietraszkiewicz-book04,Eremeyev06,Birsan-Neff-MMS-2014,Birsan-Neff-L54-2014}
\begin{equation}\label{e55}
    \begin{array}{crl}
     \boldsymbol E^e & := &  \boldsymbol Q_e^T\mathrm{Grad}_s\boldsymbol m - \boldsymbol a=
     \big(\boldsymbol m_{,\alpha}\cdot\boldsymbol d_i-\boldsymbol a_\alpha\cdot\boldsymbol d_i^0\,\big)\,
     \boldsymbol d_i^0\otimes \boldsymbol a^\alpha\,, \vspace{4pt}\\
    \boldsymbol K^e & := &  \boldsymbol Q_e^T \, \mathrm{axl}\big(\boldsymbol Q_{e,\alpha}\boldsymbol Q_e^T\big)\otimes \boldsymbol a^\alpha \, =
    \,\,\mathrm{axl}\big(\boldsymbol Q_e^T\boldsymbol Q_{e,\alpha}\big)\otimes \boldsymbol a^\alpha \vspace{4pt}\\
   & = &   \big(  \boldsymbol{d}_{2,\alpha}\cdot\boldsymbol{d}_3-  \boldsymbol{d}_{2,\alpha}^0\cdot\boldsymbol{d}_3^0\big)\,\boldsymbol{d}_1^0 \otimes\boldsymbol{a}^\alpha
       + \big( \boldsymbol{d}_{3,\alpha}\cdot\boldsymbol{d}_1-   \boldsymbol{d}_{3,\alpha}^0\cdot\boldsymbol{d}_1^0\big)\,\boldsymbol{d}_2^0 \otimes\boldsymbol{a}^\alpha
       \vspace{4pt}\\
    &  &  +\big(  \boldsymbol{d}_{1,\alpha}\cdot\boldsymbol{d}_2-   \boldsymbol{d}_{1,\alpha}^0\cdot\boldsymbol{d}_2^0\big)\,\boldsymbol{d}_3^0 \otimes\boldsymbol{a}^\alpha\,.
    \end{array}
\end{equation}
Then
\begin{equation*}
    \begin{array}{l}
     \boldsymbol c \boldsymbol K^e = \big(\boldsymbol d_1^0\otimes \boldsymbol d_2^0- \boldsymbol d_2^0\otimes \boldsymbol d_1^0\big)\boldsymbol K^e =
     \big( \boldsymbol{d}_{3,\alpha}\cdot\boldsymbol{d}_1-   \boldsymbol{d}_{3,\alpha}^0\cdot\boldsymbol{d}_1^0\big)\,\boldsymbol{d}_1^0 \otimes\boldsymbol{a}^\alpha -  \big(  \boldsymbol{d}_{2,\alpha}\cdot\boldsymbol{d}_3-  \boldsymbol{d}_{2,\alpha}^0\cdot\boldsymbol{d}_3^0\big)\,\boldsymbol{d}_2^0 \otimes\boldsymbol{a}^\alpha
     \vspace{4pt}\\
      \qquad\,\, =  \big( \boldsymbol{d}_{1}^0\otimes\boldsymbol{d}_1+   \boldsymbol{d}_{2}^0\otimes\boldsymbol{d}_2\big)\,\big(\boldsymbol{d}_{3,\alpha} \otimes\boldsymbol{a}^\alpha\big) -
      \big( \boldsymbol{d}_{1}^0\otimes\boldsymbol{d}_1^0+   \boldsymbol{d}_{2}^0\otimes\boldsymbol{d}_2^0\big)\,\big(\boldsymbol{d}_{3,\alpha}^0 \otimes\boldsymbol{a}^\alpha\big)
      \vspace{4pt}\\
       \qquad\,\, =   \boldsymbol Q_e^T\mathrm{Grad}_s\boldsymbol d_3-  \boldsymbol a  \mathrm{Grad}_s\boldsymbol n_0 \,\,= \,\,    \boldsymbol Q_e^T\mathrm{Grad}_s\boldsymbol d_3   + \boldsymbol b\,,
    \end{array}
\end{equation*}
which means
\begin{equation}\label{e56}
     \boldsymbol Q_e^T\mathrm{Grad}_s\boldsymbol d_3= \boldsymbol c \boldsymbol K^e  -  \boldsymbol b\,.
\end{equation}
We also find
\begin{equation}\label{e57}
    \begin{array}{rl}
      \boldsymbol Q_e^T\,\mathrm{Grad}_s\big(\varrho_m \boldsymbol d_3 \big) & =
      \boldsymbol Q_e^T \big[\big(\varrho_m \boldsymbol d_3\big)_{,\alpha}\otimes  \boldsymbol{a}^\alpha \big]  =
      \boldsymbol Q_e^T  \big(\varrho_m \boldsymbol d_{3,\alpha}\otimes  \boldsymbol{a}^\alpha + \varrho_{m,\alpha} \boldsymbol d_{3}\otimes  \boldsymbol{a}^\alpha \big) \vspace{4pt}\\
       & =
     \varrho_m \boldsymbol Q_e^T\mathrm{Grad}_s\boldsymbol d_3 + \varrho_{m,\alpha} \boldsymbol d_{3}^0\otimes  \boldsymbol{a}^\alpha  \,\,= \,\,
     \varrho_m (\boldsymbol c \boldsymbol K^e  -  \boldsymbol b) + \varrho_{m,\alpha} \boldsymbol n_0\otimes  \boldsymbol{a}^\alpha
    \end{array}
\end{equation}
and, analogously,
\begin{equation}\label{e58}
    \boldsymbol Q_e^T\,\mathrm{Grad}_s\big(\varrho_b \boldsymbol d_3 \big) =
     \varrho_b (\boldsymbol c \boldsymbol K^e  -  \boldsymbol b) + \varrho_{b,\alpha}\, \boldsymbol n_0\otimes  \boldsymbol{a}^\alpha\,.
\end{equation}
Using the relations \eqref{e55}--\eqref{e58} and the decomposition
$$\id_3= \boldsymbol a+ \boldsymbol n_0\otimes\boldsymbol n_0 =
 \dfrac{1}{b(x_3)}\,\big[ \boldsymbol a-x_3\,2H\boldsymbol a+x_3^2\,K\boldsymbol a   \big]
+ \boldsymbol n_0\otimes\boldsymbol n_0
$$
we can express the tensors \eqref{e53} and \eqref{e54} in terms of the shell  strain measures $\boldsymbol E^e$ and $\boldsymbol K^e\,$ by
\begin{align}\label{e59}
      \boldsymbol{E}_s = & \dfrac{1}{b(x_3)}\,\Big\{  \boldsymbol E^e+x_3\Big[ \boldsymbol E^e \big(\boldsymbol b - 2H\, \boldsymbol a\big)+\boldsymbol b  +\varrho_m(\boldsymbol c \boldsymbol K^e  -  \boldsymbol b) + \varrho_{m,\alpha} \,\boldsymbol n_0\otimes  \boldsymbol{a}^\alpha \Big]\notag
      \vspace{4pt}\\
      & \qquad \ \ \ \  + x_3^2 \Big[ \varrho_m \boldsymbol c \boldsymbol K^e \big(\boldsymbol b - 2H\, \boldsymbol a\big) + \frac{1}{2} \, \varrho_b\, (\boldsymbol c \boldsymbol K^e  -  \boldsymbol b) + \varrho_{m,\alpha}\, \boldsymbol n_0\otimes  \boldsymbol{a}^\alpha \big(\boldsymbol b - 2H\, \boldsymbol a\big)+\frac12 \,  \varrho_{b,\alpha}\,\boldsymbol n_0\otimes  \boldsymbol{a}^\alpha   \Big]
      \vspace{4pt}\\
      & \qquad \ \  \  \ +\dfrac{x_3^3}{2}\,\Big[ \varrho_b\, \boldsymbol c \boldsymbol K^e \big(\boldsymbol b - 2H\, \boldsymbol a\big) +\varrho_b\,K\,\boldsymbol a+ \varrho_{b,\alpha}\, \boldsymbol n_0\otimes  \boldsymbol{a}^\alpha \big(\boldsymbol b - 2H\, \boldsymbol a\big)
      \Big]\Big\} + [(\varrho_m-1)+x_3\varrho_b]\boldsymbol n_0\otimes  \boldsymbol n_0\,,
       \vspace{6pt}\notag\\
      \boldsymbol{\Gamma}_s = & \dfrac{1}{b(x_3)}\,\boldsymbol K^e\,\Big[ \boldsymbol a+ x_3 \big(\boldsymbol b - 2H\, \boldsymbol a\big) \Big].\notag
\end{align}

\medskip

In what follows, we shall determine the coefficients $\varrho_m$ and $\varrho_b$ by imposing the conditions \eqref{e46} on the strain tensor $\overline{\boldsymbol E} =\boldsymbol E_s\,$. In view of \eqref{e41}$_1$ and \eqref{e59}$_1$,
\begin{equation}\label{e60}
    \boldsymbol{S}_2\big(  \boldsymbol E_s \big)= 2\mu\,\boldsymbol E_s \,  +  \,2(\mu_c-\mu) \,  \mathrm{skew}\, \boldsymbol E_s \, \,   + \, \lambda \,\big( \mathrm{tr}\,\boldsymbol E_s \big)\,\id_3
\end{equation}
as well as
\begin{equation}\label{e61}
    \begin{array}{rl}
      \Big[\boldsymbol{S}_2\big(  \boldsymbol E_s \big)\Big]_{x_3=0} & = \boldsymbol{S}_2\Big(  {\boldsymbol E_s}_{\big| x_3=0}\Big),
       \vspace{6pt}\\
      \Big[ \,\dfrac{\partial}{\partial x_3}\,\boldsymbol{S}_2\big( \boldsymbol E_s \big)\,\Big]_{x_3=0} & =  2\mu\,\Big[\dfrac{\partial\boldsymbol E_s}{\partial x_3}\,\Big]_{x_3=0}  +  2(\mu_c-\mu) \,  \mathrm{skew}\, \Big[\dfrac{\partial\boldsymbol E_s}{\partial x_3}\,\Big]_{x_3=0}    +  \lambda \,\Big[\dfrac{\partial \big( \mathrm{tr}\,\boldsymbol E_s \big) }{\partial x_3}\,\Big]_{x_3=0}\id_3\,.
    \end{array}
\end{equation}
Moreover, from \eqref{e59}$_1$ we obtain
\begin{equation}\label{e62}
    \begin{array}{rl}
      \big[   \boldsymbol E_s \big]_{x_3=0} & = \boldsymbol E^e + (\varrho_m-1)\boldsymbol n_0\otimes\boldsymbol n_0\,,
       \vspace{6pt}\\
      \Big[\dfrac{\partial\boldsymbol E_s}{\partial x_3}\,\Big]_{x_3=0}  & =
       \big(\boldsymbol E^e + \boldsymbol a\big)\boldsymbol b  +\varrho_m(\boldsymbol c \boldsymbol K^e  -  \boldsymbol b) + \varrho_{m,\alpha} \,\boldsymbol n_0\otimes  \boldsymbol{a}^\alpha  + \varrho_b\,\boldsymbol n_0\otimes  \boldsymbol n_0
    \end{array}
\end{equation}
and
\begin{equation}\label{e63}
    \begin{array}{rl}
      \mathrm{tr}\,\boldsymbol E_s  =  & \dfrac{1}{b(x_3)}\,\Big\{ \mathrm{tr}\, \boldsymbol E^e+x_3\Big[ \mathrm{tr}\big(\boldsymbol E^e (\boldsymbol b - 2H\, \boldsymbol a)\big)+ 2H(1-\varrho_m) +\varrho_m\,\mathrm{tr}\,(\boldsymbol c \boldsymbol K^e)\Big]
      \vspace{4pt}\\
      & \qquad \ \  \  \ + x_3^2 \Big[ \varrho_m \,\mathrm{tr}\, \big(\boldsymbol c \boldsymbol K^e (\boldsymbol b - 2H\, \boldsymbol a)\big) + \frac{1}{2} \, \varrho_b\,\big( \mathrm{tr}\, (\boldsymbol c \boldsymbol K^e ) -  2H\big)   \Big]
      \vspace{4pt}\\
      & \qquad \ \  \  \ +\dfrac{x_3^3}{2}\,\Big[ \varrho_b\,\mathrm{tr}\, \big( \boldsymbol c \boldsymbol K^e (\boldsymbol b - 2H\, \boldsymbol a)\big) +2\varrho_b\,K   \Big]\Big\} + \big[(\varrho_m-1)+x_3\varrho_b\,\big]\,,
       \vspace{6pt}\\
       \Big[\dfrac{\partial \big( \mathrm{tr}\,\boldsymbol E_s \big) }{\partial x_3}\,\Big]_{x_3=0}  = &
        \mathrm{tr}\big(\boldsymbol E^e \boldsymbol b\big)  +2H(1-\varrho_m) +\varrho_m\,\mathrm{tr}\,(\boldsymbol c \boldsymbol K^e) +\varrho_b\,.
    \end{array}
\end{equation}
Inserting \eqref{e61}$_1$ and \eqref{e62}$_1$ in \eqref{e46}$_1$ we find
\begin{align*}
   &
       \boldsymbol n_0 \cdot \Big\{  \boldsymbol{S}_2\Big(  \big[   \boldsymbol E_s \big]_{x_3=0}\Big)\,  \cdot\boldsymbol n_0 \Big\}=0 \qquad \Leftrightarrow \vspace{4pt}\\
        & \boldsymbol n_0 \cdot \Big\{ \Big[ 2\mu\big(\boldsymbol E^e + (\varrho_m-1)\boldsymbol n_0\otimes\boldsymbol n_0  \big)  +  \,2(\mu_c-\mu) \,  \mathrm{skew}\, \boldsymbol E^e \, \,   + \, \lambda \,\big[ (\mathrm{tr}\,\boldsymbol E^e)+ (\varrho_m-1)\big]\id_3\,\Big]  \cdot\boldsymbol n_0 \Big\}=0        \qquad \Leftrightarrow \vspace{4pt}\\&
        2\mu \, (\varrho_m-1) \,+ \, \lambda \,\big[ (\mathrm{tr}\,\boldsymbol E^e)+ (\varrho_m-1)\big]=0,
     \end{align*}
which yields the expression
\begin{equation}\label{e64}
    \varrho_m=1-\dfrac{\lambda}{\lambda+2\mu}\,\big(\mathrm{tr}\,\boldsymbol E^e\big).
\end{equation}

for the coefficient $\varrho_m$. Similarly, by inserting \eqref{e62}$_2\,$, \eqref{e63}$_2$ in \eqref{e61}$_2$ and \eqref{e46}$_2$ we find
\begin{equation*}
    \begin{array}{l}
       \boldsymbol n_0 \cdot \Big\{  \Big[ \,\dfrac{\partial}{\partial x_3}\,\boldsymbol{S}_2\big( \boldsymbol E_s \big)\,\Big]_{x_3=0} \,  \cdot\boldsymbol n_0 \Big\}=0 \qquad \Leftrightarrow \vspace{4pt}\\
         \boldsymbol n_0 \cdot \Big\{ \Big[ 2\mu\big[
          \big(\boldsymbol E^e + \boldsymbol a\big)\boldsymbol b  +\varrho_m(\boldsymbol c \boldsymbol K^e  -  \boldsymbol b) + \varrho_{m,\alpha} \,\boldsymbol n_0\otimes  \boldsymbol{a}^\alpha  + \varrho_b\,\boldsymbol n_0\otimes  \boldsymbol n_0
          \big]\vspace{4pt} \\
         \qquad\quad +  \,2(\mu_c-\mu) \, \big[
           \mathrm{skew}\, \big((\boldsymbol E^e + \boldsymbol a)\boldsymbol b\big)  +\varrho_m\,\mathrm{skew}\big(\boldsymbol c \boldsymbol K^e \big) + \varrho_{m,\alpha} \,\mathrm{skew}\big(\boldsymbol n_0\otimes  \boldsymbol{a}^\alpha\big)  \, \big] \,
           \vspace{4pt}\\
           \qquad\quad + \, \lambda \,\big[  \mathrm{tr}\big(\boldsymbol E^e \boldsymbol b\big)  +2H(1-\varrho_m) +\varrho_m\,\mathrm{tr}\,(\boldsymbol c \boldsymbol K^e) +\varrho_b\,\big]\id_3\,\Big]
            \cdot\boldsymbol n_0 \Big\}=0   \qquad \Leftrightarrow
            \vspace{4pt}\\
        2\mu \, \varrho_b \,+ \, \lambda \,\big[ \, \mathrm{tr}\big(\boldsymbol E^e \boldsymbol b\big)  +2H(1-\varrho_m) +\varrho_m\,\mathrm{tr}\,(\boldsymbol c \boldsymbol K^e) +\varrho_b\,\big]=0,
     \end{array}
\end{equation*}
i.e. we can express the coefficient $\varrho_b\,$ in the form
\begin{equation}\label{e65}
    \varrho_b= -\dfrac{\lambda}{\lambda+2\mu}\,\big[  \mathrm{tr}\big(\boldsymbol E^e \boldsymbol b\big)  +\varrho_m\,\mathrm{tr}\,(\boldsymbol c \boldsymbol K^e) +2H(1-\varrho_m)\, \big].
\end{equation}
Since we consider a physically linear model, we shall neglect the quadratic terms in the deformation measures $\boldsymbol E^e$ and $  \boldsymbol K^e\,$, when replacing $\varrho_m$ given by \eqref{e64} into the relation \eqref{e65}. Thus, we obtain the following approximate  expression for $\varrho_b$ (as a simplified form of relation \eqref{e65}):
\begin{equation}\label{e66}
     \varrho_b= -\dfrac{\lambda}{\lambda+2\mu}\,\,  \mathrm{tr}\big(\,\boldsymbol E^e \boldsymbol b+ \boldsymbol c \boldsymbol K^e\,\big) .
\end{equation}
In other words, we use the approximation $\varrho_m\simeq 1$ in \eqref{e65} to obtain \eqref{e66}. Indeed, we observe that the reference values $\varrho_m^0$ and $\varrho_b^0$ of the parameters $\varrho_m$ and $\varrho_b$ are given by
\begin{equation}\label{e67}
\varrho_m^0=1,\qquad \varrho_b^0=0.
\end{equation}

\setcounter{equation}{0}

\section{Dimensionally reduced energy: analytical integration\\ through the thickness}

In order to integrate the strain energy density through the thickness, we shall use a further simplification of the form of the deformation gradient, appropriate for thin shells. Thus, we approximate the (reconstructed) deformation gradient \eqref{e52} by
\begin{equation}\label{e68}
%\begin{array}{crl}
 \widetilde{\boldsymbol F}_s  :=  \dfrac{1}{b(x_3)}\,\Big[ \mathrm{Grad}_s\boldsymbol m + x_3\mathrm{Grad}_s  \boldsymbol d_3    \Big]
  \big[ \boldsymbol a+x_3(\boldsymbol b-2H\boldsymbol a) \big] +(\varrho_m+x_3\varrho_b)\boldsymbol d_3\otimes\boldsymbol n_0\,.
%\end{array}
\end{equation}
To obtain the simplified form \eqref{e68} of \eqref{e52}, we have used the approximations $\mathrm{Grad}_s\big(\varrho_m \boldsymbol d_3 \big)\simeq \mathrm{Grad}_s  \boldsymbol d_3 \,\,$ and $\mathrm{Grad}_s\big( \varrho_b\, \boldsymbol d_3\big)\simeq \boldsymbol 0$ of the gradients as well as $\varrho_m\simeq 1$ and $\varrho_b\simeq 0$, cf. \eqref{e67}. Accordingly, for the strain tensor corresponding to $\,\boldsymbol E_s\,$ in \eqref{e59}$_1$ we find the simplified form
\begin{equation*}
    \begin{array}{crl}
  \widetilde{\boldsymbol{E}}_s \; := \;\boldsymbol Q_e^T \widetilde{\boldsymbol F}_s-\id_3 & = &
   \dfrac{1}{b(x_3)}\,\Big\{  \boldsymbol E^e+x_3\Big[ \boldsymbol E^e \big(\boldsymbol b - 2H\, \boldsymbol a\big)+ \boldsymbol c \boldsymbol K^e   \Big]
       + x_3^2 \, \boldsymbol c \boldsymbol K^e \big(\boldsymbol b - 2H\, \boldsymbol a\big) \Big\} \vspace{6pt}\\
  &   & + [(\varrho_m-1)+x_3\varrho_b]\boldsymbol n_0\otimes  \boldsymbol n_0\,,
\end{array}
\end{equation*}
which can be written as a product in the form
\begin{equation}\label{e69}
    \begin{array}{crl}
  \widetilde{\boldsymbol{E}}_s & = &
   \dfrac{1}{b(x_3)}\,\Big\{  \Big[\boldsymbol E^e + (\varrho_m-1) \boldsymbol n_0\otimes  \boldsymbol n_0\Big] +x_3\Big[ \boldsymbol E^e \big(\boldsymbol b - 2H\, \boldsymbol a\big)+ \boldsymbol c \boldsymbol K^e + \Big( 2H(1-\varrho_m)+ \varrho_b \Big) \boldsymbol n_0\otimes  \boldsymbol n_0 \Big]  \vspace{4pt}\\
  &   &\qquad \ \  \  \  
       + x_3^2 \,\Big[ \boldsymbol c \boldsymbol K^e \big(\boldsymbol b - 2H \boldsymbol a\big) + \Big( K(\varrho_m-1)- 2H\varrho_b \Big) \boldsymbol n_0\otimes  \boldsymbol n_0   \Big]  + x_3^3\, K\,\varrho_b \,\boldsymbol n_0\otimes  \boldsymbol n_0\Big\}\,.
\end{array}
\end{equation}
We denote the coefficients of $\,\boldsymbol n_0\otimes  \boldsymbol n_0\,$ appearing in \eqref{e69} by
\begin{equation}\label{e70}
    \begin{array}{l}
      A_1:= 2H(1-\varrho_m)+ \varrho_b = -\dfrac{\lambda}{\lambda+2\mu}\,\mathrm{tr}\,\Big[  \boldsymbol E^e \big(\boldsymbol b - 2H\, \boldsymbol a\big)+ \boldsymbol c \boldsymbol K^e       \Big],
      \vspace{4pt}\\
      A_2:=  K(\varrho_m-1)- 2H\varrho_b = \dfrac{\lambda}{\lambda+2\mu}\,\Big\{  \, 2H\, \mathrm{tr}\,\big[ \boldsymbol E^e\boldsymbol b  + \boldsymbol c \boldsymbol K^e \big] - K   \,\mathrm{tr}\, \boldsymbol E^e     \Big\},
    \end{array}
\end{equation}
%in view of \eqref{e64}, \eqref{e66}. 
with $\varrho_m$ and $\varrho_b$ given by \eqref{e64} and \eqref{e66}, respectively.
In the following, we want to find the expression of the strain energy density 
$$ W =W_{\mathrm{mp}}(\widetilde{\boldsymbol{E}}_s)+ W_{\mathrm{curv}}(  \boldsymbol \Gamma_s)$$ 
and to  integrate it over the thickness, according to \eqref{e20}. To this aim, we introduce the bilinear forms
\begin{equation}\label{e71}
    \begin{array}{crl}
      W_{\mathrm{mp3}}(\boldsymbol S,\boldsymbol T)& := &  \mu\,\langle\, \mathrm{sym}\, \boldsymbol S\,,\,\mathrm{sym}\, \boldsymbol T\,\rangle +  \mu_c\langle \,\mathrm{skew}\, \boldsymbol S\,,\,\mathrm{skew}\, \boldsymbol T\,\rangle +\,\dfrac{\lambda}{2}\,\big( \mathrm{tr} \boldsymbol S\big)\,\big(\mathrm{tr} \boldsymbol T\big) \vspace{4pt}\\
      & = &   \mu\,\langle\, \mathrm{dev_3\,sym}\, \boldsymbol S\,,\,\mathrm{dev_3\,sym}\, \boldsymbol T\,\rangle +  \mu_c\langle\, \mathrm{skew} \boldsymbol S\,,\,\mathrm{skew}\, \boldsymbol T\,\rangle +\,\dfrac{\kappa}{2}\,\big( \mathrm{tr} \boldsymbol S\big)\,\big(\mathrm{tr} \boldsymbol T\big), \vspace{4pt}\\
      W_{\mathrm{curv3}}(\boldsymbol S,\boldsymbol T)& := &
      \mu\,L_c^2\,\Big[ b_1\,\langle\, \mathrm{dev_3\,sym}\, \boldsymbol S\,,\,\mathrm{dev_3\,sym}\, \boldsymbol T\,\rangle +  b_2\,\langle\, \mathrm{skew} \boldsymbol S\,,\,\mathrm{skew}\, \boldsymbol T\,\rangle +b_3\,\big( \mathrm{tr} \boldsymbol S\big)\,\big(\mathrm{tr} \boldsymbol T\big) \Big],
    \end{array}
\end{equation}
for any second order tensors $\, \boldsymbol S\,,\,  \boldsymbol T\,$ in the Euclidean 3-space. We also denote the corresponding quadratic forms by
\begin{equation}\label{e72}
    W_{\mathrm{mp}}(\boldsymbol S):= W_{\mathrm{mp3}}(\boldsymbol S,\boldsymbol S),\qquad
    W_{\mathrm{curv}}(\boldsymbol S ):= W_{\mathrm{curv3}}(\boldsymbol S,\boldsymbol S),
\end{equation}
in accordance with the notation in \eqref{e16} and \eqref{e18}. Thus, using \eqref{e69} and the notation from \eqref{e70}--\eqref{e72}, we find
\begin{align}\label{e73}
   W_{\mathrm{mp}}(\widetilde{\boldsymbol{E}}_s)  =  \dfrac{1}{b^2(x_3)}\,\,W_{\mathrm{mp}} \Big( & \big[\boldsymbol E^e + (\varrho_m-1) \boldsymbol n_0\otimes  \boldsymbol n_0\big] +x_3\big[ \big(\boldsymbol E^e\boldsymbol b + \boldsymbol c \boldsymbol K^e \big) -   2H\, \boldsymbol E^e+ A_1 \boldsymbol n_0\otimes  \boldsymbol n_0 \big]  \vspace{4pt}\notag\\
  &  
       + x_3^2 \,\big[ \boldsymbol c \boldsymbol K^e\boldsymbol b   - 2H\, \boldsymbol  c \boldsymbol K^e + A_2 \boldsymbol n_0\otimes  \boldsymbol n_0   \big]  + x_3^3\, K\varrho_b \,\boldsymbol n_0\otimes  \boldsymbol n_0\Big).
   \end{align}
In order to perform the integration over the thickness, we write the right-hand side of \eqref{e73} as a polynomial in $x_3$ with the coefficients $C_k$\,, i.e.
\begin{equation}\label{e74}
    W_{\mathrm{mp}}(\widetilde{\boldsymbol{E}}_s)\,\,=\,\, \dfrac{1}{b^2(x_3)}\,  \Big( \, \sum_{k=0}^6\,C_k(x_1,x_2)\,x_3^k\,     \Big),\qquad 
\end{equation}
where
\begin{align}\label{e75}
      C_0(x_1,x_2) &:= W_{\mathrm{mp}} \Big(   \boldsymbol E^e + (\varrho_m-1) \boldsymbol n_0\otimes  \boldsymbol n_0\Big), \vspace{4pt}\notag\\
      C_1(x_1,x_2) &:=  2 \, W_{\mathrm{mp3}}   \Big(\,\boldsymbol E^e + (\varrho_m-1) \boldsymbol n_0\otimes  \boldsymbol n_0\,,\,  \big(\boldsymbol E^e\boldsymbol b + \boldsymbol c \boldsymbol K^e \big) -   2H\, \boldsymbol E^e+ A_1 \boldsymbol n_0\otimes  \boldsymbol n_0 \Big),  \vspace{4pt}\notag\\
      C_2(x_1,x_2) & :=   W_{\mathrm{mp}} \Big(  \big(\boldsymbol E^e\boldsymbol b + \boldsymbol c \boldsymbol K^e \big) -   2H\, \boldsymbol E^e+ A_1 \boldsymbol n_0\otimes  \boldsymbol n_0\Big)  \notag\\
  &   
   \qquad +  2 \, W_{\mathrm{mp3}}   \Big(\,\boldsymbol E^e + (\varrho_m-1) \boldsymbol n_0\otimes  \boldsymbol n_0\,\,,\, \boldsymbol c \boldsymbol K^e\boldsymbol b   - 2H\, \boldsymbol  c \boldsymbol K^e + A_2 \boldsymbol n_0\otimes  \boldsymbol n_0 \Big),
      \\
      C_3(x_1,x_2) & :=   2\,W_{\mathrm{mp3}} \Big(  \boldsymbol E^e + (\varrho_m-1) \boldsymbol n_0\otimes  \boldsymbol n_0\,\,,\,   K\varrho_b \,\boldsymbol n_0\otimes  \boldsymbol n_0  \Big)  \notag\\
  &   
   \qquad  +  2 \, W_{\mathrm{mp3}}   \Big(\,\big(\boldsymbol E^e\boldsymbol b + \boldsymbol c \boldsymbol K^e \big) -   2H\, \boldsymbol E^e+ A_1 \boldsymbol n_0\otimes  \boldsymbol n_0\,\,,\, \boldsymbol c \boldsymbol K^e\boldsymbol b   - 2H\, \boldsymbol  c \boldsymbol K^e + A_2 \boldsymbol n_0\otimes  \boldsymbol n_0 \Big),\notag\\
      C_4(x_1,x_2)  & :=  W_{\mathrm{mp}} \Big( \boldsymbol c \boldsymbol K^e\boldsymbol b   - 2H\, \boldsymbol  c \boldsymbol K^e + A_2 \boldsymbol n_0\otimes  \boldsymbol n_0 \Big) \notag\\
  &   
 \qquad  +  2 \, W_{\mathrm{mp3}}   \Big(\,\big(\boldsymbol E^e\boldsymbol b + \boldsymbol c \boldsymbol K^e \big) -   2H\, \boldsymbol E^e+ A_1 \boldsymbol n_0\otimes  \boldsymbol n_0\,\,,\,  K\varrho_b \,\boldsymbol n_0\otimes  \boldsymbol n_0 \Big),\vspace{4pt}\notag\\
      C_5(x_1,x_2) & :=   2 \, W_{\mathrm{mp3}}   \Big( \boldsymbol c \boldsymbol K^e\boldsymbol b   - 2H\, \boldsymbol  c \boldsymbol K^e + A_2 \boldsymbol n_0\otimes  \boldsymbol n_0\,\,,\,  K\varrho_b \,\boldsymbol n_0\otimes  \boldsymbol n_0 \Big),\vspace{4pt}\notag\\
      C_6(x_1,x_2) & :=   \, W_{\mathrm{mp}}   \Big(  K\varrho_b \,\boldsymbol n_0\otimes  \boldsymbol n_0 \Big).\notag
\end{align}
Since $x_3 \in \big(-\frac{h}{2}\,,\,\frac{h}{2}\,\big)$ for small $\,h\,$ is small, we employ the series expansion
\begin{align}\label{e76}
   \dfrac{1}{b(x_3)} \, = \,  \dfrac{1}{1-2H\,x_3+K\,x_3^2} =  1&+ 2H\,x_3+ (4H^2-K)\,x_3^2+
   (8H^3-4HK)\,x_3^3  \vspace{4pt}\notag\\
   &  +(K^2-12H^2K+16H^4)\,x_3^4+O(x_3^5)
\end{align}
and relation \eqref{e30} to compute
\begin{equation}\label{e77}
\begin{array}{c}
  \dd \int_{\Omega_h}   \,W_{\mathrm{mp}}(\widetilde{\boldsymbol{E}}_s) \,\,\mathrm{det} \big[  \nabla_x  \boldsymbol\Theta(\boldsymbol x) \big]\, \mathrm d\boldsymbol x =
    \int_{\Omega_h}   \,\Big( \, \sum_{k=0}^6\,C_k(x_1,x_2)\,x_3^k\,     \Big)\,\Big[ 1+ 2H\,x_3+ (4H^2-K)\,x_3^2    \vspace{4pt}\\
   +    (8H^3-4HK)x_3^3  +(K^2-12H^2K+16H^4)\,x_3^4+O(x_3^5)
    \Big] \,a(x_1,x_2)\,\mathrm dx_1 \mathrm d x_2 \mathrm d x_3  \vspace{4pt}\\
    = \dd\int_\omega \,\Big\{   h\,C_0+\,\dfrac{h^3}{12}\,\Big[ (4H^2-K)C_0+ 2H\,C_1+C_2 \Big]+ \,\dfrac{h^5}{80}\,\Big[ (K^2-12H^2K+16H^4)C_0    \vspace{4pt}\\
   \qquad \qquad + (8H^3-4HK)\,C_1+  (4H^2-K)\,C_2+ 2H\,C_3+C_4 \Big]
    \Big\}\,\, a(x_1,x_2) \,\mathrm dx_1 \mathrm d x_2 +O(h^7).
     \end{array}
\end{equation}
In view of \eqref{e77}, we need to find appropriate expressions for the coefficients $\,C_0\,,\,C_1\,,\,C_2\,,\,C_3\,,\,C_4\,$ defined by \eqref{e75}. We therefore denote by $\,W_{\mathrm{mixt}}(\boldsymbol S,\boldsymbol T)\,$ the bilinear form
\begin{align}\label{e78}
     W_{\mathrm{mixt}}(\boldsymbol S,\boldsymbol T) :=&   \mu\,\langle\, \mathrm{sym}\, \boldsymbol S\,,\,\mathrm{sym}\, \boldsymbol T\,\rangle +  \mu_c\langle \,\mathrm{skew}\, \boldsymbol S\,,\,\mathrm{skew}\, \boldsymbol T\,\rangle +\,\dfrac{\lambda\,\mu}{\lambda+2\mu}\,\big( \mathrm{tr} \boldsymbol S\big)\,\big(\mathrm{tr} \boldsymbol T\big) 
     \vspace{4pt}\notag\\
     \qquad\quad
       =&    \mu\,\langle\, \mathrm{dev_3\,sym}\, \boldsymbol S\,,\,\mathrm{dev_3\,sym}\, \boldsymbol T\,\rangle +  \mu_c\langle\, \mathrm{skew} \boldsymbol S\,,\,\mathrm{skew}\, \boldsymbol T\,\rangle +\,\dfrac{2\mu(2\lambda+\mu)}{3(\lambda+2\mu)}\,\big( \mathrm{tr} \boldsymbol S\big)\,\big(\mathrm{tr} \boldsymbol T\big),
       \vspace{4pt}\\
     W_{\mathrm{m}}(\boldsymbol S)  := & W_{\mathrm{mixt}}(\boldsymbol S,\boldsymbol S).\notag
    \end{align}
Observe that, since $  \,\dfrac{\kappa}{2 }\, - \,\dfrac{\lambda^2}{2(\lambda+2\mu)}\,= \,\dfrac{2\mu(2\lambda+\mu)}{3(\lambda+2\mu)}\,$,
\begin{equation}\label{e79}
    W_{\mathrm{mixt}}(\boldsymbol S,\boldsymbol T)= W_{\mathrm{mp3}}(\boldsymbol S,\boldsymbol T) - \,\dfrac{\lambda^2}{2(\lambda+2\mu)}\,\big( \mathrm{tr} \boldsymbol S\big)\,\big(\mathrm{tr} \boldsymbol T\big).
\end{equation}
Using the notation from \eqref{e71}, \eqref{e78}, we also obtain the useful relations (cf.~Appendix \ref{section:appendix})
\begin{align}\label{e80}
    &  W_{\mathrm{mp3}}\big(\boldsymbol S+\alpha\,\boldsymbol n_0\otimes  \boldsymbol n_0\,\,,\,\boldsymbol T+\beta\,\boldsymbol n_0\otimes  \boldsymbol n_0\big)= W_{\mathrm{mp3}}(\boldsymbol S,\boldsymbol T)+ \dfrac{\lambda}{2}\,\big( \alpha\,\mathrm{tr}\,\boldsymbol T+\beta\, \mathrm{tr}\,\boldsymbol S\big)+\dfrac{\lambda+2\mu}{2}\,\,\alpha\,\beta,
      \vspace{4pt}\notag\\
    &   W_{\mathrm{mp3}}\Big(\boldsymbol S-\,\dfrac{\lambda}{\lambda+2\mu}\,\big( \mathrm{tr}\,\boldsymbol S\big)\boldsymbol n_0\otimes  \boldsymbol n_0\,\,,\,\boldsymbol T+\beta\,\boldsymbol n_0\otimes  \boldsymbol n_0\Big)= W_{\mathrm{mixt}}(\boldsymbol S,\boldsymbol T),
\end{align}
which hold for all tensors of the form $\, \boldsymbol S=S_{i\gamma}\,\boldsymbol d_i^0\otimes \boldsymbol a^\gamma$ , $\, \boldsymbol T=T_{i\gamma}\,\boldsymbol d_i^0\otimes \boldsymbol a^\gamma\,$ and any coefficients $\,\alpha,\beta\in\R$. By virtue of \eqref{e80}$_2$ and \eqref{e64}, \eqref{e70}, we find
\begin{align}\label{e81}
    & W_{\mathrm{mp3}}\Big(\boldsymbol E^e + \big( \varrho_m-1\big)\boldsymbol n_0\otimes  \boldsymbol n_0\,\,,\,\boldsymbol T+\beta\,\boldsymbol n_0\otimes  \boldsymbol n_0\Big)= W_{\mathrm{mixt}}(\boldsymbol E^e,\boldsymbol T),\vspace{4pt}\notag\\
   &  W_{\mathrm{mp3}}\Big(\,\big(\boldsymbol E^e\boldsymbol b + \boldsymbol c \boldsymbol K^e \big) -   2H\, \boldsymbol E^e+ A_1 \boldsymbol n_0\otimes  \boldsymbol n_0\,\,,\,\boldsymbol T+\beta\,\boldsymbol n_0\otimes  \boldsymbol n_0\Big)= W_{\mathrm{mixt}}\big(\boldsymbol E^e\boldsymbol b + \boldsymbol c \boldsymbol K^e  -   2H\, \boldsymbol E^e\,,\,\boldsymbol T\big)\,,\vspace{4pt}\\
   &   W_{\mathrm{mp}}\Big(  \boldsymbol c \boldsymbol K^e\boldsymbol b   - 2H\, \boldsymbol  c \boldsymbol K^e + A_2 \boldsymbol n_0\otimes  \boldsymbol n_0 \Big) = W_{\mathrm{m}}\big(  \boldsymbol c \boldsymbol K^e\boldsymbol b   - 2H\, \boldsymbol c \boldsymbol K^e  \big) + \dfrac{\lambda^2}{2(\lambda+2\mu)}\Big[   \mathrm{tr}\big((\boldsymbol E^e\boldsymbol b + \boldsymbol c \boldsymbol K^e) \boldsymbol b \big) \Big]^2\notag
     \end{align}
for any tensor $\, \boldsymbol T=T_{i\gamma}\,\boldsymbol d_i^0\otimes \boldsymbol a^\gamma\,$ and any coefficient $\,\beta\in\R$.
From \eqref{e75} and \eqref{e81} we obtain
\begin{equation}\label{e82}
    \begin{array}{crl}
      C_0 & = &  W_{\mathrm{m}}\big(  \boldsymbol E^e\big), \vspace{4pt}\\
      C_1 & = & 2 \, W_{\mathrm{mixt}}  \big(\boldsymbol E^e\,,\,\boldsymbol E^e\boldsymbol b + \boldsymbol c \boldsymbol K^e  -   2H\, \boldsymbol E^e\big) = -4H\,W_{\mathrm{m}}\big(  \boldsymbol E^e\big)+ 2 \, W_{\mathrm{mixt}}  \big(\boldsymbol E^e\,,\,\boldsymbol E^e\boldsymbol b + \boldsymbol c \boldsymbol K^e \big),\vspace{4pt}\\
      C_2 & = &   W_{\mathrm{m}}  \big( \boldsymbol E^e\boldsymbol b + \boldsymbol c \boldsymbol K^e  -   2H\, \boldsymbol E^e\big)+ 2 \, W_{\mathrm{mixt}}  \big(\boldsymbol E^e\,,\,\boldsymbol c \boldsymbol K^e\boldsymbol b  -   2H\,  \boldsymbol c \boldsymbol K^e \big), \vspace{4pt}\\
      C_3 & = & 2 \, W_{\mathrm{mixt}}  \big(\boldsymbol E^e\boldsymbol b + \boldsymbol c \boldsymbol K^e  -   2H\, \boldsymbol E^e\,,\,\boldsymbol c \boldsymbol K^e\boldsymbol b  -   2H\,  \boldsymbol c \boldsymbol K^e \big), \vspace{4pt}\\
      C_4 & = &  W_{\mathrm{m}}  \big( \boldsymbol c \boldsymbol K^e\boldsymbol b  -   2H\,  \boldsymbol c \boldsymbol K^e \big)+\,\dfrac{\lambda^2}{2(\lambda+2\mu)}\,\Big[   \mathrm{tr}\big((\boldsymbol E^e\boldsymbol b + \boldsymbol c \boldsymbol K^e) \boldsymbol b \big) \Big]^2.
    \end{array}
\end{equation}
With the help of the relations \eqref{e82} we can express the two square brackets appearing in the right-hand side of \eqref{e77} by
\begin{equation*}
(4H^2-K)C_0+ 2H\,C_1+C_2 =  -K\, W_{\mathrm{m}}\big(  \boldsymbol E^e\big) + W_{\mathrm{m}}  \big( \boldsymbol E^e\boldsymbol b + \boldsymbol c \boldsymbol K^e\big) + 2 \, W_{\mathrm{mixt}}  \big(\boldsymbol E^e\,,\,\boldsymbol c \boldsymbol K^e\boldsymbol b  -   2H\,  \boldsymbol c \boldsymbol K^e \big)
\end{equation*}
and
\begin{equation}\label{e83}
    \begin{array}{l}
      (K^2-12H^2K+16H^4)C_0
    + (8H^3-4HK)\,C_1+  (4H^2-K)\,C_2+ 2H\,C_3+C_4  = \vspace{4pt}\\
     \qquad\quad =
     -K\,  W_{\mathrm{m}}  \big( \boldsymbol E^e\boldsymbol b +  \boldsymbol c \boldsymbol K^e\big) +  W_{\mathrm{m}} \big((\boldsymbol E^e\boldsymbol b + \boldsymbol c \boldsymbol K^e) \boldsymbol b\, \big) + \,\dfrac{\lambda^2}{2(\lambda+2\mu)}\,\Big[   \mathrm{tr}\big((\boldsymbol E^e\boldsymbol b + \boldsymbol c \boldsymbol K^e) \boldsymbol b \big) \Big]^2  \vspace{4pt}\\
     \qquad\quad =
     -K\,  W_{\mathrm{m}}  \big( \boldsymbol E^e\boldsymbol b +  \boldsymbol c \boldsymbol K^e\big) +  W_{\mathrm{mp}} \big((\boldsymbol E^e\boldsymbol b + \boldsymbol c \boldsymbol K^e) \boldsymbol b \,\big),
    \end{array}
\end{equation}
respectively. Inserting \eqref{e83} into \eqref{e77} (and neglecting the terms of order $O(h^7)$) we obtain the following result for the integration:
\begin{equation}\label{e84}
   \begin{array}{l}
  \dd \int_{\Omega_h}   \,W_{\mathrm{mp}}(\widetilde{\boldsymbol{E}}_s) \,\,\mathrm{det} \big[  \nabla_x  \boldsymbol\Theta(\boldsymbol x) \big]\, \mathrm d\boldsymbol x  \vspace{4pt}\\
  \qquad\qquad =
    \dd\int_{\omega}   \,\, \Big[  \Big(h-K\,\dfrac{h^3}{12}\Big)\,
    W_{\mathrm{m}}\big(  \boldsymbol E^e\big)    +  \Big(\dfrac{h^3}{12}\,-K\,\dfrac{h^5}{80}\Big)\,
     W_{\mathrm{m}}  \big( \boldsymbol E^e\boldsymbol b +  \boldsymbol c \boldsymbol K^e\big) \vspace{4pt}\\
   \qquad\qquad\quad\ \ \;  \;\;\, +
     \dfrac{h^3}{12}\,\,2\, W_{\mathrm{mixt}}  \big(\boldsymbol E^e\,,\,\boldsymbol c \boldsymbol K^e\boldsymbol b  -   2H\,  \boldsymbol c \boldsymbol K^e \big)+ \,\dfrac{h^5}{80}\,\,
     W_{\mathrm{mp}} \big((\boldsymbol E^e\boldsymbol b + \boldsymbol c \boldsymbol K^e) \boldsymbol b \big)
      \Big] a(x_1,x_2)      \,\mathrm dx_1 \mathrm d x_2  .
     \end{array}
\end{equation}

%\bigskip

Analogously, we integrate the curvature part of the strain energy density, according to \eqref{e20}, \eqref{e30}, \eqref{e59}$_2$ and \eqref{e76}:
\begin{equation}\label{e85}
    \begin{array}{l}
       \dd \int_{\Omega_h}   \,W_{\mathrm{curv}}( \boldsymbol{\Gamma}) \,\,\mathrm{det} \big[  \nabla_x  \boldsymbol\Theta(\boldsymbol x) \big]\, \mathrm d\boldsymbol x =  \int_{\Omega_h}  \,\,W_{\mathrm{curv}}\big(\boldsymbol K^e+ x_3 \,(  \boldsymbol K^e\boldsymbol b   - 2H\,  \boldsymbol K^e)\,\big)
         \,\dfrac{a(x_1,x_2)}{b(x_3)}\,\,\mathrm d\boldsymbol x 
         \vspace{4pt}\\
    \qquad =
     \dd\int_{\Omega_h}   \,\Big( D_0+ D_1\,x_3 +D_2\,x_3^2\,     \Big)\,\Big[ 1+ 2H\,x_3+ (4H^2-K)\,x_3^2
   +    (8H^3-4HK)x_3^3 
   \vspace{4pt}\\
   \qquad\qquad\quad\;\; +(K^2-12H^2K+16H^4)\,x_3^4+O(x_3^5)
    \Big] \,a(x_1,x_2)\,\mathrm dx_1 \,\mathrm d x_2\, \mathrm d x_3  \vspace{4pt}\\
     \qquad     = \dd\int_\omega \,\Big\{   h\,D_0+\,\dfrac{h^3}{12}\,\Big[ (4H^2-K)D_0+ 2H\,D_1+D_2 \Big]+ \,\dfrac{h^5}{80}\,\Big[ (K^2-12H^2K+16H^4)\,D_0    
     \vspace{4pt}\\
    \qquad\qquad\quad\; + (8H^3-4HK)\,D_1+  (4H^2-K)\,D_2  \Big]
    \Big\}\,\, a(x_1,x_2) \,\mathrm dx_1 \mathrm d x_2 +O(h^7),
    \end{array}
\end{equation}
where we have denoted by $D_k$ the coefficients of $x_3^k$ ($k=0,1,2$) in the expression
\begin{equation}\label{e86}
    \begin{array}{c}
      W_{\mathrm{curv}}\big(\boldsymbol K^e+ x_3 \,(  \boldsymbol K^e\boldsymbol b   - 2H\,  \boldsymbol K^e)\,\big)=  D_0(x_1,x_2)+ D_1(x_1,x_2)\,x_3 +D_2(x_1,x_2)\,x_3^2\,,\qquad\text{i.e.}
      \vspace{4pt}\\
       D_0 := W_{\mathrm{curv}}\big(  \boldsymbol K^e\big), \quad
      D_1 := 2 \, W_{\mathrm{curv3}} \big(\boldsymbol K^e\,, \, \boldsymbol K^e\boldsymbol b   - 2H\,  \boldsymbol K^e\big),\quad
      D_2 :=   W_{\mathrm{curv}} \big(  \boldsymbol K^e\boldsymbol b   - 2H\,  \boldsymbol K^e\big).
    \end{array}
\end{equation}
We write the coefficients of $\,\dfrac{h^3}{12}\,$ and $\,\dfrac{h^5}{80}\,$ in \eqref{e85} with the help of \eqref{e86}$_{2,3,4}$ as
\begin{equation}\label{e87}
     \begin{array}{l}
     (4H^2-K)D_0+ 2H\,D_1+D_2 = -K\, W_{\mathrm{curv}}\big(\boldsymbol K^e \big) + W_{\mathrm{curv}}\big(\boldsymbol K^e \boldsymbol b \,  \big) \qquad\text{and}
     \vspace{4pt}\\
     (K^2-12H^2K+16H^4)\,D_0   + (8H^3-4HK)\,D_1+  (4H^2-K)\,D_2 
     \vspace{4pt}\\
    \qquad\qquad = K^2\, W_{\mathrm{curv}}\big(\boldsymbol K^e \big)
        -4HK\,  W_{\mathrm{curv3}}\big(\boldsymbol K^e\,,\,\boldsymbol K^e \boldsymbol b \,  \big)+ (4H^2-K)\, W_{\mathrm{curv}}\big(\boldsymbol K^e \boldsymbol b \,  \big)
        \vspace{4pt}\\
    \qquad\qquad
    = -K\, W_{\mathrm{curv}}\big(\boldsymbol K^e\boldsymbol b \,   \big) + W_{\mathrm{curv}}\big(\boldsymbol K^e \boldsymbol b^2  \big).
    \end{array}
\end{equation}
Inserting \eqref{e87} into \eqref{e85} (and neglecting the terms of order $O(h^7)$) we arrive at the following result for the integration of the curvature part:
\begin{equation}\label{e88}
     \begin{array}{l}
       \dd \int_{\Omega_h}   \,W_{\mathrm{curv}}( \boldsymbol{\Gamma}) \,\,\mathrm{det} \big[  \nabla_x  \boldsymbol\Theta(\boldsymbol x) \big]\, \mathrm d\boldsymbol x \vspace{4pt}\\
  \qquad =
    \dd\int_{\omega}   \,\, \Big[  \Big(h-K\,\dfrac{h^3}{12}\Big)\,
   W_{\mathrm{curv}}\big(\boldsymbol K^e \big)    +  \Big(\dfrac{h^3}{12}\,-K\,\dfrac{h^5}{80}\Big)\,
      W_{\mathrm{curv}}\big(\boldsymbol K^e \boldsymbol b \,  \big)  + \,\dfrac{h^5}{80}\,\,
      W_{\mathrm{curv}}\big(\boldsymbol K^e \boldsymbol b^2  \big)
      \Big] \,a      \,\mathrm dx_1 \mathrm d x_2  .
      \end{array}
\end{equation}
%where the quadratic form  is defined by \eqref{e72}. \bigskip

\bigskip
\noindent Finally, from \eqref{e20}, \eqref{e84} and \eqref{e88} we see that the total energy functional $I$ for shells has the form
\begin{equation}\label{e89}
    I= \int_{\omega}   \,\, \Big[  \,
   W_{\mathrm{memb}}\big(\boldsymbol E^e\,,\,\boldsymbol K^e \big)    +
      W_{\mathrm{bend}}\big(\boldsymbol K^e    \big)
      \Big] \,a(x_1,x_2)      \,\mathrm dx_1 \mathrm d x_2\,,
\end{equation}
where the  membrane part $\,W_{\mathrm{memb}}\big(\boldsymbol E^e\,,\,\boldsymbol K^e \big) \,$ and the bending--curvature part $\,W_{\mathrm{bend}}\big(\boldsymbol K^e    \big)\,$ of the shell energy density are given by
\begin{align}\label{e90}
      W_{\mathrm{memb}}\big(\boldsymbol E^e\,,\,\boldsymbol K^e \big)&=  \Big(h-K\,\dfrac{h^3}{12}\Big)\,
    W_{\mathrm{m}}\big(  \boldsymbol E^e\big)    +  \Big(\dfrac{h^3}{12}\,-K\,\dfrac{h^5}{80}\Big)\,
     W_{\mathrm{m}}  \big( \boldsymbol E^e\boldsymbol b +  \boldsymbol c \boldsymbol K^e\big) 
     \vspace{4pt}\notag\\
  & \qquad +
     \dfrac{h^3}{12}\,\,2\, W_{\mathrm{mixt}}  \big(\boldsymbol E^e\,,\,\boldsymbol c \boldsymbol K^e\boldsymbol b  -   2H\,  \boldsymbol c \boldsymbol K^e \big)+ \,\dfrac{h^5}{80}\,\,
     W_{\mathrm{mp}} \big((\boldsymbol E^e\boldsymbol b + \boldsymbol c \boldsymbol K^e) \boldsymbol b \,\big), 
     \vspace{4pt}\\
      W_{\mathrm{bend}}\big(\boldsymbol K^e    \big) &=  \Big(h-K\,\dfrac{h^3}{12}\Big)\,
   W_{\mathrm{curv}}\big(\boldsymbol K^e \big)    +  \Big(\dfrac{h^3}{12}\,-K\,\dfrac{h^5}{80}\Big)\,
      W_{\mathrm{curv}}\big(\boldsymbol K^e \boldsymbol b \,  \big)  + \,\dfrac{h^5}{80}\,\,
      W_{\mathrm{curv}}\big(\boldsymbol K^e \boldsymbol b^2  \big),\medskip\notag
\end{align}
with $\, W_{\mathrm{m}}\,,\, W_{\mathrm{mixt}}\,,\, W_{\mathrm{mp}}\,$ and  $\,W_{\mathrm{curv}}\,$ given by \eqref{e72} and \eqref{e78}.

\section{Conclusions and final remarks}

Starting with a three-dimensional Cosserat model, we perform a  dimensional reduction and derive a two-dimensional shell model. Beginning with the 8-parameter ansatz \eqref{e51}, we determine subsequently two of the parameters (namely $ \varrho_m $ in \eqref{e64} and $ \varrho_b $ in \eqref{e66}), using the generalized plane stress condition \eqref{e42}. Thus, we finally obtain a 6-parameter model for Cosserat shells.

The dimensionally reduced strain energy density for shells is obtained by explicit analytic integration through the thickness, retaining all the terms up to the order $ O(h^5) $. This energy density for shells admits the additive split \eqref{e89} into the membrane part $\,W_{\mathrm{memb}}\big(\boldsymbol E^e\,,\,\boldsymbol K^e \big) \,$ and the bending--curvature part $\,W_{\mathrm{bend}}\big(\boldsymbol K^e    \big)\,$ and is expressed as a quadratic function of the elastic shell strain tensor $ \,\boldsymbol E^e\, $ and the elastic shell bending-curvature tensor $ \,\boldsymbol K^e\, $. These strain measures are commonly used in  the general nonlinear shell theory, see e.g. \cite{Libai98,Pietraszkiewicz04,Pietraszkiewicz-book04}. The coefficients of this strain energy density for shells also depend on the initial curvature of the reference midsurface via the second fundamental tensor $ \boldsymbol b $, the mean curvature $ H=\frac{1}{2}\mathrm{tr}\,\boldsymbol b\,  $ and the Gau{\ss}  curvature $ K=\mathrm{det}\,\boldsymbol b \,$, see \eqref{e90}.

We remark that the specific form \eqref{e90} of the strain energy density satisfies all invariance properties required by the local symmetry group of isotropic shells, established in the general 6-parameter theory of shells in \cite[Section 9]{Eremeyev06}.

To show existence results for the nonlinear Cosserat shell model established here, one might apply the results from \cite{Birsan-Neff-MMS-2014}. Thus, this model is indeed viable for applications and can directly be implemented to numerically solve nonlinear shell problems with large rotations and curved initial configurations, cf.~\cite{Sander-Neff-Birsan-16}.

\appendix
\section{Appendix}
\label{section:appendix}

In this appendix we prove the relations \eqref{e80} by a straightforward calculation: in view of the definition \eqref{e71},
\begin{align}
      W_{\mathrm{mp3}}\big(\boldsymbol S+\alpha\,\boldsymbol n_0\otimes  \boldsymbol n_0\,\,,\,\boldsymbol T+\beta\,\boldsymbol n_0\otimes  \boldsymbol n_0\big)&=  \mu\,\langle\, \mathrm{sym}\, \boldsymbol S+\alpha\,\boldsymbol n_0\otimes  \boldsymbol n_0\,,\,\mathrm{sym}\, \boldsymbol T+\beta\,\boldsymbol n_0\otimes  \boldsymbol n_0\,\rangle
      \vspace{4pt}\\
        &\quad +  \mu_c\langle \,\mathrm{skew}\, \boldsymbol S\,,\,\mathrm{skew}\, \boldsymbol T\,\rangle +\,\dfrac{\lambda}{2}\,\big( \mathrm{tr} \boldsymbol S+\alpha\big)\,\big(\mathrm{tr}\, \boldsymbol T+\beta\big).\notag
    \end{align}

Since $\, \boldsymbol S=S_{i\gamma}\,\boldsymbol d_i^0\otimes \boldsymbol a^\gamma$ , $\, \boldsymbol T=T_{i\gamma}\,\boldsymbol d_i^0\otimes \boldsymbol a^\gamma\,$, we find $\,
\langle\, \mathrm{sym}\, \boldsymbol S \,,\, \boldsymbol n_0\otimes  \boldsymbol n_0\,\rangle=0\,$ , $\,
\langle\, \mathrm{sym}\, \boldsymbol T \,,\, \boldsymbol n_0\otimes  \boldsymbol n_0\,\rangle=0\,$ and thus
\begin{align}
      W_{\mathrm{mp3}}\big(\boldsymbol S+\alpha\,\boldsymbol n_0\otimes  \boldsymbol n_0\,\,,\,\boldsymbol T+\beta\,\boldsymbol n_0\otimes  \boldsymbol n_0\big)&=  \mu\,\langle\, \mathrm{sym}\, \boldsymbol S\,,\,\mathrm{sym}\, \boldsymbol T\,\rangle +\,\dfrac{\lambda}{2}\,\big( \mathrm{tr} \boldsymbol S\big)\,\big(\mathrm{tr} \boldsymbol T\big)  +
      \mu_c\langle \,\mathrm{skew}\, \boldsymbol S\,,\,\mathrm{skew}\, \boldsymbol T\,\rangle
       \vspace{4pt}\notag\\
      & \qquad  + \mu\,\alpha\,\beta\, \langle\, \boldsymbol n_0\otimes  \boldsymbol n_0 \,,\, \boldsymbol n_0\otimes  \boldsymbol n_0\,\rangle+\dfrac{\lambda}{2}\,\big( \alpha\,\mathrm{tr}\,\boldsymbol T+\beta\, \mathrm{tr}\,\boldsymbol S\big)+\dfrac{\lambda}{2}\,\,\alpha\,\beta 
       \vspace{4pt}\notag\\
    &  = W_{\mathrm{mp3}}(\boldsymbol S,\boldsymbol T)+ \dfrac{\lambda}{2}\,\big( \alpha\,\mathrm{tr}\,\boldsymbol T+\beta\, \mathrm{tr}\,\boldsymbol S\big)+\dfrac{\lambda+2\mu}{2}\,\,\alpha\,\beta\,,
    \end{align}
which means that the relation \eqref{e80}$_1$ holds true, for any $\, \boldsymbol S\,,\,  \boldsymbol T\,$  and $\,\alpha,\beta\in\R$.

If we write \eqref{e80}$_1$ with $\,\alpha=\,\dfrac{-\lambda}{\lambda+2\mu}\,\big( \mathrm{tr}\,\boldsymbol S\big)\,$ , then using formula \eqref{e79} we obtain
\begin{equation} 
    \begin{array}{l}
       W_{\mathrm{mp3}}\Big(\boldsymbol S-\,\dfrac{\lambda}{\lambda+2\mu}\,\big( \mathrm{tr}\,\boldsymbol S\big)\boldsymbol n_0\otimes  \boldsymbol n_0\,\,,\,\boldsymbol T+\beta\,\boldsymbol n_0\otimes  \boldsymbol n_0\Big) \vspace{4pt}\\
      \qquad\qquad =
       W_{\mathrm{mp3}}(\boldsymbol S,\boldsymbol T)
       + \dfrac{\lambda}{2}\,\Big( \,\dfrac{-\lambda}{\lambda+2\mu}\,\big( \mathrm{tr}\,\boldsymbol S\big)\big(\mathrm{tr}\,\boldsymbol T\big)+\beta\, \mathrm{tr}\,\boldsymbol S\Big) +\dfrac{\lambda+2\mu}{2}\cdot\dfrac{-\lambda}{\lambda+2\mu}\,\big( \mathrm{tr}\,\boldsymbol S\big)\,\beta
         \vspace{4pt}\\
        \qquad\qquad =W_{\mathrm{mp3}}(\boldsymbol S,\boldsymbol T) - \,\dfrac{\lambda^2}{2(\lambda+2\mu)}\,\big( \mathrm{tr} \boldsymbol S\big)\,\big(\mathrm{tr} \boldsymbol T\big)=
        W_{\mathrm{mixt}}(\boldsymbol S,\boldsymbol T),
    \end{array}
\end{equation}
which shows the relation \eqref{e80}$_2$.

\bigskip\bigskip
\noindent
\small{\textbf{Acknowledgements}\quad
	This research has been funded by the Deutsche Forschungsgemeinschaft (DFG, German Research Foundation) -- Project no. 415894848 (M. B\^{\i}rsan and P. Neff). The work of I.D. Ghiba has been supported by a grant of the Romanian National Authority for Scientific Research and Innovation, CNCS-UEFISCDI, project number PN-III-P1-1.1-TE- 2016-2314.

\bibliographystyle{plain} %plain

\bibliography{literatur1,literatur_Birsan}

\begin{thebibliography}{10}

\bibitem{Birsan-Neff-MMS-2014}
M.~B\^{\i}rsan and P.~Neff.
\newblock Existence of minimizers in the geometrically non-linear 6-parameter
  resultant shell theory with drilling rotations.
\newblock {\em Math. Mech. Solids}, 19(4):376--397, 2014.

\bibitem{Birsan-Neff-L54-2014}
M.~B\^{\i}rsan and P.~Neff.
\newblock Shells without drilling rotations: {A} representation theorem in the
  framework of the geometrically nonlinear 6-parameter resultant shell theory.
\newblock {\em Int. J. Engng. Sci.}, 80:32--42, 2014.

\bibitem{Birsan-Neff-L57-2016}
M.~B\^{\i}rsan and P.~Neff.
\newblock On the dislocation density tensor in the {C}osserat theory of elastic
  shells.
\newblock In K.~Naumenko and M.~Assmus, editors, {\em Advanced Methods of
  Continuum Mechanics for Materials and Structures}, Advanced Structured
  Materials 60, pages 391--413. Springer Science+Business Media, Singapore,
  2016.

\bibitem{Birsan-Neff-L58-2017}
M.~B\^{\i}rsan and P.~Neff.
\newblock Analysis of the deformation of {C}osserat elastic shells using the
  dislocation density tensor.
\newblock In F.~dell'Isola~et al., editor, {\em Advanced Methods of Continuum
  Mechanics for Materials and Structures}, Advanced Structured Materials 69,
  pages 13--30. Springer Nature, Singapore, 2017.

\bibitem{Ramm97}
M.~Bischoff and E.~Ramm.
\newblock Shear deformable shell elements for large strains and rotations.
\newblock {\em Int. J. Num. Meth. Engrg.}, 40:4427--4449, 1997.

\bibitem{Ramm00}
M.~Bischoff and E.~Ramm.
\newblock On the physical significance of higher order kinematic and static
  variables in a three-dimensional shell formulation.
\newblock {\em Int. J. Solids Struct.}, 37:6933--6960, 2000.

\bibitem{Ramm94}
M.~Braun, M.~Bischoff, and E.~Ramm.
\newblock Nonlinear shell formulations for complete three-dimensional
  constitutive laws including composites and laminates.
\newblock {\em Comp. Mech.}, 15:1--18, 1994.

\bibitem{Ramm92}
N.~B\"uchter and E.~Ramm.
\newblock Shell theory versus degeneration-a comparison in large rotation
  finite element analysis.
\newblock {\em Int. J. Num. Meth. Engrg.}, 34:39--59, 1992.

\bibitem{Chernykh80}
K.~Chernykh.
\newblock Nonlinear theory of isotropically elastic thin shells.
\newblock {\em Mechanics of Solids, Transl. of Mekh. Tverdogo Tela},
  15(2):118--127, 1980.

\bibitem{Pietraszkiewicz-book04}
J.~Chr\'o\'scielewski, J.~Makowski, and W.~Pietraszkiewicz.
\newblock {\em Statics and Dynamics of Multifold Shells: Nonlinear Theory and
  Finite Element Method (in Polish).}
\newblock Wydawnictwo IPPT PAN, Warsaw, 2004.

\bibitem{Ciarlet00}
P.G. Ciarlet.
\newblock {\em Mathematical {E}lasticity, {V}ol. {III}: {T}heory of {S}hells.}
\newblock North-Holland, Amsterdam, first edition, 2000.

\bibitem{Pietraszkiewicz04}
V.A. Eremeyev and W.~Pietraszkiewicz.
\newblock The nonlinear theory of elastic shells with phase transitions.
\newblock {\em J. Elasticity}, 74:67--86, 2004.

\bibitem{Eremeyev06}
V.A. Eremeyev and W.~Pietraszkiewicz.
\newblock Local symmetry group in the general theory of elastic shells.
\newblock {\em J. Elasticity}, 85:125--152, 2006.

\bibitem{Libai98}
A.~Libai and J.G. Simmonds.
\newblock {\em The {N}onlinear {T}heory of {E}lastic {S}hells.}
\newblock Cambridge University Press, Cambridge, 1998.

\bibitem{Neff_plate04_cmt}
P.~Neff.
\newblock A geometrically exact {C}osserat-shell model including size effects,
  avoiding degeneracy in the thin shell limit. {P}art {I}: {F}ormal dimensional
  reduction for elastic plates and existence of minimizers for positive
  {C}osserat couple modulus.
\newblock {\em Cont. Mech. Thermodynamics}, 16(6 (DOI
  10.1007/s00161-004-0182-4)):577--628, 2004.

\bibitem{Neff_Habil04}
P.~Neff.
\newblock {\em Geometrically exact {C}osserat theory for bulk behaviour and
  thin structures. {M}odelling and mathematical analysis.}
\newblock Signatur HS 7/0973. Habilitationsschrift, Universit\"ats- und
  Landesbibliothek, Technische Universit\"at Darmstadt, Darmstadt, 2004.

\bibitem{Neff_zamm06}
P.~Neff.
\newblock The {C}osserat couple modulus for continuous solids is zero viz the
  linearized {C}auchy-stress tensor is symmetric.
\newblock {\em Z. Angew. Math. Mech.}, 86:892--912, 2006.

\bibitem{Neff_Edinb06}
P.~Neff.
\newblock Existence of minimizers for a finite-strain micromorphic elastic
  solid.
\newblock {\em Proc. Roy. Soc. Edinb.}, 136A:997--1012, 2006.

\bibitem{Neff_plate07_m3as}
P.~Neff.
\newblock A geometrically exact planar {C}osserat shell-model with
  microstructure: Existence of minimizers for zero {C}osserat couple modulus.
\newblock {\em Math. Mod. Meth. Appl. Sci.}, 17:363--392, 2007.

\bibitem{Birsan-Neff-Ost_L56-2015}
P.~Neff, M.~B\^{\i}rsan, and F.~Osterbrink.
\newblock Existence theorem for geometrically nonlinear {C}osserat micropolar
  model under uniform convexity requirements.
\newblock {\em J. Elasticity}, 121:119--141, 2015.

\bibitem{Neff_Chelminski_ifb07}
P.~Neff and K.~Che{\l}mi\'nski.
\newblock A geometrically exact {C}osserat shell-model for defective elastic
  crystals. {J}ustification via {$\Gamma$}-convergence.
\newblock {\em Interfaces and Free Boundaries}, 9:455--492, 2007.

\bibitem{Neff_Hong_Reissner08}
P.~Neff, K.-I. Hong, and J.~Jeong.
\newblock The {R}eissner-{M}indlin plate is the {$\Gamma$}-limit of {C}osserat
  elasticity.
\newblock {\em Math. Mod. Meth. Appl. Sci.}, 20:1553--1590, 2010.

\bibitem{Neff_curl08}
P.~Neff and I.~M\"unch.
\newblock Curl bounds {Grad} on {${\rm SO}(3)$}.
\newblock {\em ESAIM: Control, Optimisation and Calculus of Variations},
  14:148--159, 2008.

\bibitem{Neff-Fischle-14}
P.~Neff, Y.~Nakatsukasa, and A.~Fischle.
\newblock A logarithmic minimization property of the unitary polar factor in
  the spectral and {F}robenius norms.
\newblock {\em SIAM J. Matrix Anal. Appl.}, 25:1132--1154, 2014.

\bibitem{Pietraszkiewicz85}
W.~Pietraszkiewicz.
\newblock {\em Finite {R}otations in {S}tructural {M}echanics.}
\newblock Number~19 in Lectures Notes in Engineering. Springer, Berlin, 1985.

\bibitem{Pietraszkiewicz09}
W.~Pietraszkiewicz and V.A. Eremeyev.
\newblock On natural strain measures of the non-linear micropolar continuum.
\newblock {\em Int. J. Solids Struct.}, 46:774--787, 2009.

\bibitem{Pietraszkiewicz14}
W.~Pietraszkiewicz and V.~Konopi\'nska.
\newblock Drilling couples and refined constitutive equations in the resultant
  geometrically non-linear theory of elastic shells.
\newblock {\em Int. J. Solids Struct.}, 51:2133--2143, 2014.

\bibitem{Reissner74}
E.~Reissner.
\newblock Linear and nonlinear theory of shells.
\newblock In Y.C. Fung and E.E. Sechler, editors, {\em Thin Shell Structures.},
  pages 29--44. Prentice-Hall, Englewood Cliffs, New Jersey, 1974.

\bibitem{Ramm96}
D.~Roehl and E.~Ramm.
\newblock Large elasto-plastic finite element analysis of solids and shells
  with the enhanced assumed strain concept.
\newblock {\em Int. J. Solids Struct.}, 33:3215--3237, 1996.

\bibitem{Sander-Neff-Birsan-16}
O.~Sander, P.~Neff, and M.~B\^{\i}rsan.
\newblock Numerical treatment of a geometrically nonlinear planar {C}osserat
  shell model.
\newblock {\em Computational Mechanics}, 57:817--841, 2016.

\bibitem{Sansour98c}
C.~Sansour and J.~Bocko.
\newblock On hybrid stress, hybrid strain and enhanced strain finite element
  formulations for a geometrically exact shell theory with drilling degrees of
  freedom.
\newblock {\em Int. J. Num. Meth. Engrg.}, 43:175--192, 1998.

\bibitem{Schmidt85}
R.~Schmidt.
\newblock Polar decomposition and finite rotation vector in first order finite
  elastic strain shell theory.
\newblock In W.~Pietraszkiewicz, editor, {\em Finite {R}otations in
  {S}tructural {M}echanics}, number~19 in Lecture Notes in Engineering.
  Springer, Berlin, 1985.

\bibitem{Steigmann12}
D.J. Steigmann.
\newblock Extension of {K}oiter's linear shell theory to materials exhibiting
  arbitrary symmetry.
\newblock {\em Int. J. Engng. Sci.}, 51:216--232, 2012.

\bibitem{Steigmann13}
D.J. Steigmann.
\newblock Koiter's shell theory from the perspective of three-dimensional
  nonlinear elasticity.
\newblock {\em J. Elasticity}, 111:91--107, 2013.

\bibitem{Tambaca-19}
J.~Tamba\'{c}a, M.~Ljulj, and Z.~Tutek.
\newblock A {N}aghdi type nonlinear model for shells with little regularity.
\newblock {\em paper communicated at GAMM Annual Meeting}, Vienna
  (Austria):February 18--22, 2019.

\bibitem{Tambaca-16}
J.~Tamba\'{c}a and Z.~Tutek.
\newblock A new linear {N}aghdi type shell model for shells with little
  regularity.
\newblock {\em Applied Mathematical Modelling}, 40:10549--10562, 2016.

\bibitem{Zhilin06}
P.A. Zhilin.
\newblock {\em Applied Mechanics -- Foundations of Shell Theory (in Russian)}.
\newblock State Polytechnical University Publisher, Sankt Petersburg, 2006.

\end{thebibliography}

\end{document}